\documentclass[article]{IEEEtran}
\usepackage{cite}
\usepackage{color}
\usepackage[pdftex]{graphicx}
\graphicspath{{fig/}{jpeg/}}
\usepackage[cmex10]{amsmath}
\usepackage{amssymb}
\usepackage{algorithm}
\usepackage{algorithmic}
\usepackage{amsmath,amssymb,lipsum}
\usepackage{hyperref}
\usepackage{comment}
\usepackage{graphicx}
\usepackage[font=small]{caption}
\usepackage{subcaption}
\usepackage{setspace}

\input{mysymbol.sty}
\usepackage{needspace}

\newcommand{\myparagraph}[1]{\needspace{1\baselineskip}\medskip\noindent {\it #1.}}

\renewcommand{\QED}{\hfill\ensuremath{\blacksquare}}
\usepackage{theorem}
\newtheorem{thm}{Theorem}
\newtheorem{lemma}{Lemma}
\newtheorem{proposition}{Proposition}
\theoremstyle{definition}
\newtheorem{assumption}{Assumption}

\title{{Network Newton--Part I: Algorithm and Convergence}}

\author{Aryan Mokhtari, Qing Ling and Alejandro Ribeiro 
\thanks{{Work in this paper is supported by NSF CAREER CCF-0952867, ONR N00014-12-1-0997, and NSFC 61004137. Aryan Mokhtari and Alejandro Ribeiro are with the Department of Electrical and Systems Engineering, University of Pennsylvania, 200 South 33rd Street, Philadelphia, PA 19104, USA. Email: \{aryanm, aribeiro\}@seas.upenn.edu. Qing Ling is with the Department of Automation, University of Science and Technology of China, 96 Jinzhao Road, Hefei, Anhui, 230026, China. Email: qingling@mail.ustc.edu.cn. Part of the results in this paper appeared in \cite{NNAsilomar} and \cite{NNICASSP}. {This paper expands the results and presents convergence proofs that are referenced in \cite{NNAsilomar} and \cite{NNICASSP}.} 
}}
}

\begin{document}

\maketitle
\thispagestyle{empty}

\begin{abstract}

We study the problem of minimizing a sum of convex objective functions where the components of the objective are available at different nodes of a network and nodes are allowed to only communicate with their neighbors. The use of distributed gradient methods is a common approach to solve this problem. Their popularity notwithstanding, these methods exhibit slow convergence and a consequent large number of communications between nodes to approach the optimal argument because they rely on first order information only. This paper proposes the network Newton (NN) method as a distributed algorithm that incorporates second order information. This is done via distributed implementation of approximations of a suitably chosen Newton step. The approximations are obtained by truncation of the Newton step's Taylor expansion. This leads to a family of methods defined by the number $K$ of Taylor series terms kept in the approximation. When keeping $K$ terms of the Taylor series, the method is called NN-$K$ and can be implemented through the aggregation of information in $K$-hop neighborhoods. Convergence to a point close to the optimal argument at a rate that is at least linear is proven and the existence of a tradeoff between convergence time and the distance to the optimal argument is shown. Convergence rate, several practical implementation matters, and numerical analyses are presented in a {companion paper \cite{NN-part2}.} 

\end{abstract}



\section{Introduction}\label{sec_Introduction}

Distributed optimization algorithms are used to solve the problem of minimizing a global cost function over a set of nodes in situations where the objective function is defined as a sum of local functions. To be more precise, consider a variable $\bbx\in\reals^p$ and a connected network containing $n$ agents each of which has access to a local function $f_i:\reals^p\to\reals$. The agents cooperate in minimizing the aggregate cost function $f:\reals^p\to\reals$ taking values $f(\bbx) := \sum_{i=1}^{n}f_i(\bbx)$. I.e., agents cooperate in solving the global optimization problem
\begin{equation}\label{original_optimization_problem1}
  \bbx^* \ :=\ \argmin_{\bbx}  f(\bbx) 
      \  =\ \argmin_{\bbx}\sum_{i=1}^{n} f_i(\bbx).
\end{equation}
Problems of this form arise often in, e.g., decentralized control systems \cite{Bullo2009,Cao2013-TII,LopesEtal8}, wireless systems \cite{Ribeiro10,Ribeiro12,rabbat2004decentralized}, sensor networks \cite{Schizas2008-1,KhanEtal10,cRabbatNowak04}, and large scale machine learning \cite{bekkerman2011scaling,Tsianos2012-allerton-consensus,Cevher2014}. {In the latter case, distributed formulations are efficient in dealing with very large datasets where it is desirable to split training sets into smaller subsamples that are assigned to different servers \cite{low2012distributed}.} In this paper we assume that the local costs $f_i$ are twice differentiable and strongly convex. Therefore, the aggregate cost function $f$ is also twice differentiable and strongly convex.

There are different algorithms to solve \eqref{original_optimization_problem1} in a distributed manner. The most popular choices are decentralized gradient descent (DGD) \cite{Nedic2009,Jakovetic2014-1,YuanQing, Shi2014}, distributed implementations of the alternating direction method of multipliers \cite{Schizas2008-1,cQingRibeiroADMM14,BoydEtalADMM11,Shi2014-ADMM}, and decentralized dual averaging (DDA) \cite{Duchi2012,cTsianosEtal12}. Although there are substantial differences between them, these methods can be generically abstracted as combinations of local descent steps followed by variable exchanges and averaging of information among neighbors. A feature common to all of these algorithms is the slow convergence rate in ill-conditioned problems since they operate on first order information only. This is not surprising because gradient descent methods in centralized settings where the aggregate function gradient is available at a single server have the same difficulties in problems with skewed curvature {[see Chapter 9 of \cite{boyd}.]}

This issue is addressed in centralized optimization by Newton's method that uses second order information to determine a descent direction adapted to the objective's curvature [see Chapter 9 of \cite{boyd}]. In general, second order methods are not available in distributed settings because distributed approximations of Newton steps are difficult to devise. In the particular case of flow optimization problems, these approximations are possible when operating in the dual domain \cite{ZarghamEtal14,ZarghamEtal14-2,WeiEtal13,cZarghamEtal13-2}. As would be expected, these methods result in large reductions of convergence times.

Our goal here is to develop approximate Newton's methods to solve \eqref{original_optimization_problem1} in distributed settings where agents have access to their local functions only and exchange variables with neighboring agents. We do so by introducing Network Newton (NN), a method that relies on distributed approximations of Newton steps for the global cost function $f$ to accelerate convergence of DGD. We begin the paper with an alternative formulation of \eqref{original_optimization_problem1} and a brief discussion of DGD (Section \ref{sec:problem}). We then introduce a reinterpretation of DGD as an algorithm that utilizes gradient descent to solve a penalized version of \eqref{original_optimization_problem1} in lieu of the original optimization problem (Section \ref{sec_penalty_interpretation}). This reinterpretation explains convergence of DGD to a neighborhood of the optimal solution. The volume of this neighborhood is given by the relative weight of the penalty function and the original objective which is controlled by a penalty coefficient. 

If gradient descent on the penalized function finds an approximate solution to the original problem, the same solution can be found with a much smaller number of iterations by using Newton's method. Alas, distributed computation of Newton steps requires global communication between all nodes in the network and is therefore impractical (Section \ref{subsec_nn}). To resolve this issue we approximate the Newton step of the penalized objective function by truncating the Taylor series expansion of the exact Newton step (Section \ref{sec_dist_appro_Newton_step}). This results in a family of methods indexed by the number of terms of the Taylor expansion that are kept in the approximation. The method that results from keeping $K$ of these terms is termed NN-$K$. A fundamental observation here is that the Hessian of the penalized function has a sparsity structure that is the same sparsity pattern of the graph. Thus, when computing terms in the Hessian inverse expansion, the first order term is as sparse as the graph, the second term is as sparse as the two hop neighborhood, and, in general, the $k$-th term is as sparse as the $k$-hop neighborhood of the graph. Thus, implementation of the NN-$K$ method requires aggregating information from $K$ hops away. Increasing $K$ makes NN-$K$ arbitrarily close to Newton's method at the cost of increasing the communication overhead of each iteration.

Convergence of NN-$K$ to the optimal argument of the penalized objective is established (Section \ref{sec:convergence_analysis}). We do so by establishing several auxiliary bounds on the eigenvalues of the matrices involved in the definition of the method (Propositions \ref{eigenvalue_bounds}-\ref{error_matrix_proposition} and Lemma \ref{Hessian_inverse_eigenvalue_bounds_lemma}). Of particular note, we show that a measure of the error between the Hessian inverse approximation utilized by NN-$K$ and the actual inverse Hessian decays exponentially with the method index $K$. This exponential decrease hints that using a small value of $K$ should suffice in practice. Convergence is formally claimed in Theorem \ref{linear_convergence} that shows the convergence rate is at least linear. It follows from this convergence analysis that larger penalty coefficients result in faster convergence that comes at the cost of increasing the distance between the optimal solutions of the original and penalized objectives. The convergence guarantees established in this paper are not better than the corresponding guarantees for DGD. These advantages are established in a companion paper where we further show that the sequence of penalized objective function values generated by NN-$K$ has a convergence rate that is quadratic in a specific interval. This quadratic phase holds for all $K$ and can be made arbitrarily large by increasing $K$ \cite{NN-part2}. Numerical results in \cite{NN-part2} establish the advantages of NN-$K$ in terms of number of iterations and communications steps relative to DGD and establish that using $K=1$ or $K=2$ tends to work best in practice.

\myparagraph{\bf Notation} Vectors are written as $\bbx\in\reals^n$ and matrices as $\bbA\in\reals^{n\times n}$. Given $n$ vectors $\bbx_i$, the vector $\bby=[\bbx_1;\ldots;\bbx_n]$ represents a stacking of the elements of each individual $\bbx_i$. The null space of matrix $\bbA$ is denoted by $\text{null}(\bbA)$ and the span of a vector by $\text{span}(\bbx)$. We use $\|\bbx\|$ to denote the Euclidean norm of vector $\bbx$ and $\|\bbA\|$ to denote the Euclidean norm of matrix $\bbA$. The gradient of a function $f(\bbx)$ is denoted as $\nabla f(\bbx)$ and the Hessian matrix is denoted as $\nabla^2 f(\bbx)$. The $i$-th largest eigenvalue of matrix $\bbA$ is denoted by $\mu_{i}(\bbA)$.


%
\section{Distributed Gradient Descent} \label{sec:problem}

The network that connects the $n$ agents is assumed connected, symmetric, and specified by the neighborhoods $\mathcal{N}_i$ that contain the list of nodes than can communicate with $i$ for $i=1,\ldots,n$. In the problem in \eqref{original_optimization_problem1} agent $i$ has access to the local objective function $f_{i}(\bbx)$ and agents cooperate to minimize the global cost $f(\bbx)$. This specification is more naturally formulated by an alternative representation of   \eqref{original_optimization_problem1} in which node $i$ selects a local decision vector $\bbx_{i}\in \reals^p$. Nodes then try to achieve the minimum of their local objective functions $f_{i}(\bbx_{i})$, while keeping their variables equal to the variables $\bbx_j$ of neighboring nodes $j\in \ccalN_i$. This alternative formulation can be written as 
\begin{align}\label{original_optimization_problem2}
   \{\bbx_i^*\}_{i=1}^n\ := \
   &\argmin_{\{\bbx_{i}\}_{i=1}^n} \ \sum_{i=1}^{n}\ f_{i}(\bbx_{i}), \nonumber\\ 
   &\text{\ s.t.}  \ \bbx_{i}=\bbx_{j}, 
                   \quad \text{for all\ } i, j\in\ccalN_i .
\end{align} 
Since the network is connected, the constraints $\bbx_{i}=\bbx_{j}$ for all $i$ and $j\in\ccalN_i$ imply that \eqref{original_optimization_problem1} and \eqref{original_optimization_problem2} are equivalent in the sense that we have $\bbx_i^*=\bbx^*$ for all $i$. This must be the case because for a connected network the constraints $\bbx_{i}=\bbx_{j}$ for all $i$ and $j\in\ccalN_i$ collapse the feasible space of \eqref{original_optimization_problem2} to a hyperplane in which all local variables are equal. When all local variables are equal, the objectives in \eqref{original_optimization_problem1} and \eqref{original_optimization_problem2} coincide and, in particular, so do their optima.

DGD is an established distributed method to solve \eqref{original_optimization_problem2} which relies on the introduction of nonnegative weights $w_{ij}\geq0$ that are not null if and only if  $j=i$ or if $j\in \mathcal{N}_i$. Letting $t\in\naturals$ be a discrete time index and $\alpha$ a given stepsize, DGD is defined by the recursion
\begin{equation}\label{gd_iteration}
   \bbx_{i,t+1} 
      = \sum_{j=1}^{n} w_{ij}\bbx_{j,t}-\alpha\nabla f_i({\bbx_{i,t}}), 
      \qquad i=1,\ldots,n.
\end{equation}
Since $w_{ij}=0$ when $j\neq i$ and $j\notin \mathcal{N}_i$, it follows from \eqref{gd_iteration} that each agent $i$ updates its estimate $\bbx_i$ of the optimal vector $\bbx^*$ by performing an average over the estimates $\bbx_{j,t}$ of its neighbors $j\in \mathcal{N}_i$ and its own estimate $\bbx_{i,t}$, and descending through the negative local gradient $-\nabla f_i(\bbx_{i,t})$. DGD is a distributed method because to implement \eqref{gd_iteration}, node $i$ exchanges variables with neighboring nodes only. 

The weights in \eqref{gd_iteration} cannot be arbitrary. To express conditions on the set of allowable weights define the matrix $\bbW\in\reals^{n\times n}$ with entries $w_{ij}$. We require the weights to be symmetric, i.e., $w_{ij}=w_{ji}$ for all $i,j$, and such that the weights of a given node sum up to 1, i.e., $\sum_{j=1}^{n} w_{ij}=1$ for all $i$. If the weights sum up to 1 we must have $\bbW\bbone=\bbone$ which implies that $\bbI-\bbW$ is rank deficient. It is also customary to require the rank of $\bbI-\bbW$ to be exactly equal to $n-1$ so that the null space of $\bbI-\bbW$ is $\text{null}(\bbI-\bbW)=\text{span}(\bbone)$. We therefore have the following three restrictions on the matrix $\bbW$,
\begin{equation}\label{eqn_conditions_on_weights}
   \bbW^T=\bbW, \quad 
   \bbW\bbone=\bbone, \quad
   \text{null}(\bbI-\bbW)=\text{span}(\bbone).
\end{equation}
If the considions in \eqref{eqn_conditions_on_weights} are true, it is possible to show that \eqref{gd_iteration} approaches the solution of \eqref{original_optimization_problem1} in the sense that $\bbx_{i,t}\approx\bbx^*$ for all $i$ and large $t$, \cite{Nedic2009}. The accepted  interpretation of why \eqref{gd_iteration} converges is that nodes are gradient descending towards their local minima because of the term $-\alpha\nabla f_i({\bbx_{i,t}})$ but also perform an average of neighboring variables $\sum_{j=1}^{n} w_{ij}\bbx_{j,t}$. This latter consensus operation drives the agents to agreement. In the following section we show that \eqref{gd_iteration} can be alternatively interpreted as a penalty method.

%
\subsection{Penalty method interpretation}\label{sec_penalty_interpretation}

It is illuminating to define matrices and vectors so as to rewrite \eqref{gd_iteration} as a single equation. To do so define the vectors $\bby := \left[\bbx_{1}; \dots ; \bbx_{n}\right]$ and $\bbh(\bby):= \left[\nabla f_{1}(\bbx_{1}); \dots ; \nabla f_{n}(\bbx_{n})\right]$. Vector $\bby\in\reals^{np}$ concatenates the local vectors $\bbx_{i}$, and the vector $\bbh(\bby)\in\reals^{np}$ concatenates the gradients of the local functions $f_i$ taken with respect to the local variable $\bbx_i$. Notice that $\bbh(\bby)$ is {\it not} the gradient of $f(\bbx)$ and that a vector $\bby$ with $\bbh(\bby)=\bbzero$ does {\it not} necessarily solve \eqref{original_optimization_problem1}. To solve \eqref{original_optimization_problem1} we need to have $\bbx_{i}=\bbx_{j}$ for all $i$ and $j$ with $\sum_{i=1}^n\nabla f_{i}(\bbx_{i})=\bbzero$. In any event, to rewrite \eqref{gd_iteration} we also define the matrix $\bbZ:= \bbW \otimes \bbI\in\reals^{np\times np}$ as the Kronecker product of the weight matrix $\bbW\in\reals^{n\times n}$ and the identity matrix $\bbI\in\reals^{p\times p}$. It is then ready to see that \eqref{gd_iteration} is equivalent to
\begin{equation}\label{new_formulation2}
{   \bby_{t+1} = \bbZ \bby_t - \alpha\bbh(\bby_t)
               = \bby_{t} - \big[  (\bbI -\bbZ) \bby_t +\alpha \bbh(\bby_t)\big],}
\end{equation}
where in the second equality we added and subtracted $\bby_t$ and regrouped terms. Inspection of \eqref{new_formulation2} reveals that the DGD update formula at step $t$ is equivalent to a (regular) gradient descent algorithm being used to solve the program
\begin{equation}\label{centralized_opt_problem}
    \bby^* := \argmin\ F(\bby) 
           := \min \frac{1}{2}\ \bby^{T}(\bbI -\bbZ)\ \bby 
              + \alpha\sum_{i=1}^{n} f_i(\bbx_i).
\end{equation}
Indeed, given the definition of the function $F(\bby):= (1/2)\bby^{T}(\bbI -\bbZ)\ \bby + \alpha\sum_{i=1}^{n} f_i(\bbx_i)$ it follows that the gradient of $F(\bby)$ at $\bby=\bby_t$ is given by
\begin{equation}\label{eqn_gradient_definition}
    \bbg_t \ :=\ \nabla F(\bby_t) 
           \  =\ (\bbI -\bbZ) \bby_t +\alpha \bbh(\bby_t).
\end{equation}
Using \eqref{eqn_gradient_definition} we rewrite \eqref{new_formulation2} as $\bby_{t+1}=\bby_{t}- \bbg_t$ and conclude that DGD descends along the negative gradient of $F(\bby)$ with unit stepsize. The expression in \eqref{gd_iteration} is just a distributed implementation of gradient descent that uses the gradient in \eqref{eqn_gradient_definition}. To confirm that this is true, observe that the $i$th element of the gradient $\bbg_{t}=[\bbg_{i,t};\ldots;\bbg_{i,t}]$ is given by
\begin{equation}\label{local_gradient}
\bbg_{i,t}=(1-w_{ii})\bbx_{i,t} - \sum_{j\in \mathcal{N}_i} w_{ij} \bbx_{j,t}+\alpha \nabla f_{i}(\bbx_{i,t}).
\end{equation}
The gradient descent iteration $\bby_{t+1}=\bby_{t}- \bbg_t$ is then equivalent to \eqref{gd_iteration} if we entrust node $i$ with the implementation of the descent $\bbx_{i,t+1}=\bbx_{i,t}- \bbg_{i,t}$, where, we recall, $\bbx_{i,t}$ and $\bbx_{i,t+1}$ are the $i$th components of the vectors $\bby_t$ and $\bby_{t+1}$. Observe that the local gradient component $\bbg_{i,t}$ can be computed using local information and the $\bbx_{j,t}$ iterates of its neighbors $j\in \ccalN_i$. This is as it should be, because the descent $\bbx_{i,t+1}=\bbx_{i,t}- \bbg_{i,t}$ is equivalent to \eqref{gd_iteration}.

Is it a good idea to descend on $F(\bby)$ to solve \eqref{original_optimization_problem1}? To some extent. Since we know that the null space of $\bbI-\bbW$ is $\text{null}(\bbI-\bbW)=\text{span}(\bbone)$ and that $\bbZ= \bbW \otimes \bbI$ we know that the span of $\bbI-\bbZ$ is $\text{null}(\bbI-\bbZ)=\text{span}(\bbone\otimes\bbI)$. Thus, we have that $(\bbI-\bbZ)\bby=\bb0$ holds if and only if $\bbx_1=\dots=\bbx_n$. Since the matrix $\bbI-\bbZ$ is positive semidefinite -- because it is stochastic and symmetric --, the same is true of the square root matrix $({\bbI-\bbZ})^{1/2}$. Therefore, we have that the optimization problem in \eqref{original_optimization_problem2} is equivalent to the optimization problem 
\begin{align}\label{original_optimization_new_notation}
   \tby^*\ :=\ &\argmin_{\bbx}\ \sum_{i=1}^{n}\ f_i(\bbx_i), \nonumber \\
               &\text{\ s.t.} \quad ({\bbI-\bbZ})^{1/2}  \bby =\bb0.
\end{align} 
Indeed, for $\bby=[\bbx_1;\ldots;\bbx_n]$ to be feasible in \eqref{original_optimization_new_notation} we must have $\bbx_1=\dots=\bbx_n$ because $\text{null}[(\bbI-\bbZ)^{1/2}]=\text{span}(\bbone\otimes\bbI)$ as already argued. This is the same constraint imposed in \eqref{original_optimization_problem2} from where it follows that we must have $\tby^*=[\bbx_1^{*};\ldots;\bbx_n^{*}]$ with $\bbx_i^{*}=\bbx^*$ for all $i$. The unconstrained minimization in \eqref{centralized_opt_problem} is a penalty version of \eqref{original_optimization_new_notation}. The penalty function associated with the constraint $({\bbI-\bbZ})^{1/2}  \bby =\bb0$ is the squared norm $(1/2)\|({\bbI-\bbZ})^{1/2}  \bby \|^2$ and the corresponding penalty coefficient is $1/\alpha$. Inasmuch as the penalty coefficient $1/\alpha$ is sufficiently large, the optimal arguments $\bby^*$ and $\tby^*$ are not too far apart. 

The reinterpretation of \eqref{gd_iteration} as a penalty method demonstrates that DGD is an algorithm that finds the optimal solution of \eqref{centralized_opt_problem}, not \eqref{original_optimization_new_notation} or its equivalent original formulations in  \eqref{original_optimization_problem1} and \eqref{original_optimization_problem2}. Using a fixed $\alpha$ the distance between  $\bby^*$ and $\tby^*$ is of order $O(\alpha)$, \cite{YuanQing}. To solve \eqref{original_optimization_new_notation} we need to introduce a rule to progressively decrease $\alpha$. In this paper we exploit the reinterpretation of \eqref{new_formulation2} as a method to minimize \eqref{centralized_opt_problem} to propose an approximate Newton algorithm that can be implemented in a distributed manner. We explain this algorithm in the following section.

%
\section{Network Newton}\label{subsec_nn}

Instead of solving \eqref{centralized_opt_problem} with a gradient descent algorithm as in DGD, we can solve \eqref{centralized_opt_problem} using Newton's method. To implement Newton's method we need to compute the Hessian $\bbH_t:=\nabla^2 F(\bby_t)$ of $F$ evaluated at $\bby_t$ so as to determine the Newton step $\bbd_t:=-\bbH_t^{-1}\bbg_t$. Start by differentiating twice in \eqref{centralized_opt_problem} in order to write $\bbH_t$ as
\begin{equation}\label{Hessian}
   \bbH_t := \nabla^2 F(\bby_t) = \bbI-\bbZ +\alpha \bbG_t,
\end{equation}
where the matrix $\bbG_t \in \reals^{np\times np}$ is a block diagonal matrix formed by blocks $\bbG_{ii,t}\in \reals^{p\times p}$ containing the Hessian of the $i$th local function,
\begin{equation}\label{G_form}
   \bbG_{ii,t} = \nabla^2 f_i(\bbx_{i,t})  .
\end{equation}
It follows from \eqref{Hessian} and \eqref{G_form} that the Hessian $\bbH_t$ is block sparse with blocks $\bbH_{ij,t}\in \reals^{p\times p}$ having the sparsity pattern of $\bbZ$, which is the sparsity pattern of the graph. {The diagonal blocks are of the form $\bbH_{ii,t}=(1-w_{ii})\bbI +  \alpha \nabla^2 f_i(\bbx_{i,t}) $ and the off diagonal blocks are not null only when $j\in \mathcal{N}_i$ in which case $\bbH_{ij,t}=w_{ij}\bbI$.} 

While the Hessian $\bbH_t$ is sparse, the inverse $\bbH_t$ is not. It is the latter that we need to compute the Newton step $\bbd_t:=\bbH_t^{-1}\bbg_t$. To overcome this problem we split the diagonal and off diagonal blocks of $\bbH_t$ and rely on a Taylor's expansion of the inverse. To be precise, write $\bbH_t=\bbD_t - \bbB$ where the matrix $\bbD_t$ is defined as 
\begin{equation}\label{diagonal_matrix}
   \bbD_t := \alpha \bbG_t + 2\ ( \bbI  -  \diag(\bbZ))
          := \alpha \bbG_t + 2\ ( \bbI  -  \bbZ_{d}),
\end{equation}
where in the second equality we defined $\bbZ_{d}:=\diag(\bbZ)$ for future reference. Since the diagonal weights must be $w_{ii}<1$, the matrix $\bbI  -  \bbZ_{d}$ is positive definite. The same is true of the block diagonal matrix $\bbG_t$ because the local functions are assumed strongly convex. Therefore, the matrix $\bbD_t$ is block diagonal and positive definite. The $i$th diagonal block $\bbD_{ii,t}\in\reals^p$ of $\bbD_t$ can be computed and stored by node $i$ as $\bbD_{ii,t}= \alpha \nabla^2 f_{i}(\bbx_{i,t}) + 2(1-w_{ii})\bbI $. To have $\bbH_t=\bbD_t - \bbB$ we must define $\bbB:=\bbD_t-\bbH_t$.  Considering the definitions of $\bbH_t$ and $\bbD_t$ in \eqref{Hessian} and \eqref{diagonal_matrix}, it follows that 
\begin{equation}\label{non_diagona_matrix}
   \bbB =  \bbI - 2\bbZ_{d} +\bbZ.
\end{equation}
Observe that $\bbB$ is independent of time and depends on the weight matrix $\bbZ$ only. As in the case of the Hessian $\bbH_t$, the matrix $\bbB$ is block sparse with blocks $\bbB_{ij}\in \reals^{p\times p}$ having the sparsity pattern of $\bbZ$, which is the sparsity pattern of the graph. Node $i$ can compute the diagonal blocks $\bbB_{ii}=(1-w_{ii})\bbI$ and the off diagonal blocks $\bbB_{ij}=w_{ij}\bbI$ using the local information about its own weights.

Proceed now to factor $\bbD_t^{1/2}$ from both sides of the splitting relationship to write $\bbH_t = \bbD_t ^{{1}/{2}} (  \bbI - \bbD_t ^{-{1}/{2}}\bbB\bbD_t ^{-{1}/{2}} ) \bbD_t^{{1}/{2}}$. When we consider the Hessian inverse $\bbH^{-1}$, we can use the Taylor series  $(\bbI-\bbX)^{-1}= \sum_{j=0}^{\infty} \bbX^{j}$ with $\bbX=\bbD_t^{-{1}/{2}}  \bbB   \bbD_t^{-{1}/{2}}$ to write  
\begin{equation}\label{exact_Hessian_inverse}
   \bbH_t^{-1} = \bbD_t^{-1/2} 
                \sum_{k=0}^{\infty} \left(\bbD_t^{-1/2}  
                  \bbB   \bbD_t^{-1/2}\right)^{k}\ \bbD_t^{-1/2}.
\end{equation}
Observe that the sum in \eqref{exact_Hessian_inverse} converges if the absolute value of all the eigenvalues of the matrix $\bbD^{-{1}/{2}}  \bbB   \bbD^{-{1}/{2}} $ are strictly less  than 1. For the time being we assume this to be the case but we will prove that this is true in Section \ref{sec:convergence_analysis}. When the series converge, we can use truncations of this series to define approximations to the Newton step as we explain in the following section.

%
\subsection{Distributed approximations of the Newton step}\label{sec_dist_appro_Newton_step} 

Network Newton (NN) is defined as a family of algorithms that rely on truncations of the series in \eqref{exact_Hessian_inverse}. The {$K$th} member of this family, NN-$K$, considers the first $K+1$ terms of the series to define the approximate Hessian inverse
\begin{equation}\label{Hessian_inverse_approximation}
   \hbH_t^{(K)^{-1}} :=     \bbD_t^{-1/2}  \  \sum_{k=0}^{K} \left(  \bbD_t^{-1/2}  \bbB   \bbD_t^{-1/2}            \right)^{k}     \ \bbD_t^{-1/2}.
\end{equation}
NN-$K$ uses the approximate Hessian $\hbH_t^{(K)^{-1}}$ as a curvature  correction matrix that is used in lieu of the exact Hessian inverse $\bbH^{-1}$ to estimate the Newton step. I.e., instead of descending along the Newton step $\bbd_t:=-\bbH_t^{-1}\bbg_t$ we descend along the NN-$K$ step $\bbd_t^{(K)}:=-\hbH_t^{(K)^{-1}}\bbg_t$, which we intend as an approximation of $\bbd_t$. Using the explicit expression for $\hbH_t^{(K)^{-1}}$ in \eqref{Hessian_inverse_approximation} we write the NN-$K$ step as
\begin{equation}\label{Hessian_approximation_iteration}
\bbd_t^{(K)} = -\  \bbD_t^{-1/2}  \  \sum_{k=0}^{K} \left(  \bbD_t^{-1/2}  \bbB   \bbD_t^{-1/2}            \right)^{k}     \ \bbD_t^{-1/2}\ \bbg_t,
\end{equation}
where, we recall, the vector $\bbg_t$ is the gradient of objective function $F(\bby)$ defined in \eqref{eqn_gradient_definition}. The NN-$K$ update formula can then be written as
\begin{equation}\label{update_formula_NN}
   \bby_{t+1}=\bby_t+\eps\  \bbd_t^{(K)}.
\end{equation}
The algorithm defined by recursive application of \eqref{update_formula_NN} can be implemented in a distributed manner because the truncated series in \eqref{Hessian_inverse_approximation} has a local structure controlled by the parameter $K$. To explain this statement better define the components $\bbd^{(K)}_{i,t}\in\reals^p$ of the NN-$K$ step $\bbd^{(K)}_{t}=[\bbd^{(K)}_{1,t};\ldots;\bbd^{(K)}_{n,t}]$. A distributed implementation of \eqref{update_formula_NN} requires that node $i$ computes $\bbd^{(K)}_{i,t}$ so as to implement the local descent $\bbx_{i,t+1}=\bbx_{i,t} + \eps\bbd^{(K)}_{i,t}$. The key observation here is that the step component $\bbd^{(K)}_{i,t}$ can indeed be computed through local operations. Specificially, begin by noting that as per the definition of the NN-$K$ descent direction in \eqref{Hessian_approximation_iteration} the sequence of NN descent directions satisfies 
\begin{equation}
   \bbd_t^{(k+1)} = \bbD_t^{-1}\bbB \bbd_t^{(k)} -\bbD_t^{-1}\bbg_t
                  = \bbD_t^{-1}\left(\bbB \bbd_t^{(k)} - \bbg_t \right).
\end{equation}
Then observe that since the matrix $\bbB$ has the sparsity pattern of the graph, this recursion can be decomposed into local components
\begin{equation}\label{local_descent}
  \bbd_{i,t}^{(k+1)} 
      = \bbD_{ii,t}^{-1}\bigg(
              \sum_{j\in \mathcal{N}_i,j=i} \bbB_{ij} \bbd_{j,t}^{(k)} 
              - \bbg_{i,t}\bigg),
\end{equation}
The matrix $\bbD_{ii,t}=\alpha \nabla^2 f_{i}(\bbx_{i,t}) + 2(1-w_{ii})\bbI$ is stored and computed at node $i$. The gradient component $\bbg_{i,t}=(1-w_{ii})\bbx_{i,t} - \sum_{j\in \mathcal{N}_i} w_{ij} \bbx_{j,t}+\alpha \nabla f_{i}(\bbx_{i,t})$ is also stored and computed at $i$. Node $i$ can also evaluate the values of the matrix blocks  $\bbB_{ii}=(1-w_{ii})\bbI$  and $\bbB_{ij}=w_{ij}\bbI$. Thus, if the NN-$k$ step components $\bbd_{j,t}^{(k)}$ are available at neighboring nodes $j$, node $i$ can then determine the  NN-$(k+1)$ step component $\bbd_{i,t}^{(k+1)}$ upon being communicated that information. 

The expression in \eqref{local_descent} represents an iterative computation embedded inside the NN-$K$  recursion in \eqref{update_formula_NN}. For each time index $t$, we compute the local component of the NN-$0$ step $\bbd_{i,t}^{(0)}=-\bbD_{ii,t}^{-1}\bbg_{i,t}$. Upon exchanging this information with neighbors we use \eqref{local_descent} to determine the NN-$1$ step components $\bbd_{i,t}^{(1)}$. These can be exchanged and plugged in \eqref{local_descent} to compute $\bbd_{i,t}^{(2)}$. Repeating this procedure $K$ times, nodes ends up having determined their NN-$K$ step component $\bbd_{i,t}^{(K)}$ .

%
\begin{algorithm}[t]{\small
\caption{Network Newton-$K$ method at node $i$}\label{algo_NN1} 
\begin{algorithmic}[1] {
\REQUIRE  Initial iterate $\bbx_{i,0}$. Weights $w_{ij}$. Penalty coefficient $\alpha$.
\STATE $\bbB$ matrix blocks: $\bbB_{ii}=(1-w_{ii})\bbI$ and $\bbB_{ij}=w_{ij}\bbI$
\FOR {$t=0,1,2,\ldots$}
   \STATE $\bbD$ matrix block: $\bbD_{ii,t}= \alpha \nabla^2 f_{i}(\bbx_{i,t}) + 2(1-w_{ii})\bbI $
   \STATE Exchange iterates $\bbx_{i,t}$ with neighbors $\displaystyle{j\in \mathcal{N}_i}$.
   \STATE Gradient:    
          $\displaystyle{
          \bbg_{i,t} = (1-w_{ii})\bbx_{i,t} 
                       - \sum_{j\in \mathcal{N}_i} w_{ij} \bbx_{j,t}
                       +\alpha \nabla f_{i}(\bbx_{i,t})}$  
   \STATE Compute NN-0 descent direction $\bbd_{i,t}^{(0)}=-\bbD_{ii,t}^{-1}\bbg_{i,t}$\\ 
   \FOR  {$k=  0, \ldots, K-1$ } 
      \STATE Exchange elements $\bbd_{i,t}^{(k)}$ of the NN-$k$ step with neighbors
      \STATE NN-$(k+1)$ step:
             $\displaystyle{  
             \bbd_{i,t}^{(k+1)} = \bbD_{ii,t}^{-1}
             			\bigg[
                                   \sum_{j\in \mathcal{N}_i,j=i}\bbB_{ij} \bbd_{j,t}^{(k)} 
                                    - \bbg_{i,t}\bigg]}$        
   \ENDFOR
          %
          \STATE Update local iterate: 
          $\displaystyle{\bbx_{i,t+1}=\bbx_{i,t} +\eps\ \bbd_{i,t}^{(K)}}$.
\ENDFOR}
\end{algorithmic}}\end{algorithm}

%
The resulting NN-$K$ method is summarized in Algorithm \ref{algo_NN1}. The descent iteration in \eqref{update_formula_NN} is implemented in Step 11. Implementation of this descent requires access to the NN-$K$ descent direction $\bbd_{i,t}^{(K)}$ which is computed by the loop in steps 6-10. Step $6$ initializes the loop by computing the NN-0 step $\bbd_{i,t}^{(0)}=-\bbD_{ii,t}^{-1}\bbg_{i,t}$. The core of the loop is in Step 9 which corresponds to the recursion in \eqref{local_descent}. Step 8 stands for the variable exchange that is necessary to implement Step 9. After $K$ iterations through this loop, the NN-$K$ descent direction $\bbd_{i,t}^{(K)}$ is computed and can be used in Step 11. Both, steps 6 and 9, require access to the local gradient component $\bbg_{i,t}$. This is evaluated in Step 5 after receiving the prerequisite information from neighbors in Step 4. Steps 1 and 3 compute the blocks $\bbB_{ii,t}$, $\bbB_{ij,t}$, and $\bbD_{ii,t}$ that are also necessary in steps 6 and 9.

\section{Convergence Analysis}\label{sec:convergence_analysis}

In this section we show that as time progresses the sequence of objective function values $F(\bby_t)$ [cf. \eqref{centralized_opt_problem}] approaches the optimal objective function value $F(\bby^{*})$. In proving this claim we make the following assumptions.

%
\begin{assumption}\label{ass_weight_bounds} There exists constants $0\leq\delta\leq\Delta<1$ that lower and upper bound the diagonal weights for all $i$, 
\begin{equation}\label{bounds_for_local_weights}
   0\leq \delta  \leq w_{ii} \leq \Delta <1,  \qquad  i=1,\ldots,n .
\end{equation}\end{assumption}

%
\begin{assumption}\label{convexity_assumption} 
The local objective functions $f_i(\bbx)$ are twice differentiable and the eigenvalues of the local objective function Hessians are bounded with positive constants $0<m\leq M<\infty$, i.e. 
\begin{equation}\label{local_hessian_eigenvlaue_bounds}
m\bbI\preceq \nabla^2 f_i(\bbx)\preceq M\bbI.
\end{equation}
\end{assumption}

%
\begin{assumption}\label{Lipschitz_assumption} The local objective function Hessians $\nabla^2 f_i(\bbx)$ are Lipschitz continuous with respect to the Euclidian norm with parameter $L$. I.e., for all $\bbx, \hbx \in \reals^p$, it holds
\begin{equation}
   \| \nabla^2 f_i(\bbx)-\nabla^2 f_i(\hbx) \| \ \leq\  L\ \| \bbx- \hbx \|.
\end{equation}
\end{assumption}

%
Notice that the lower bound in Assumption \ref{ass_weight_bounds} is more a definition than a constraint since we may have $\delta=0$. This is not recommendable as it is implies that the weight $w_{ii}$ assigned to the local variable $\bbx_i$ in \eqref{gd_iteration} is null, but nonetheless allowed. The upper bound $\Delta<1$ on the weights $w_{ii}$ is true for all connected networks as long as neighbors $j\in\ccalN_i$ are assigned nonzero weights $w_{ij}>0$. This is because the matrix $\bbW$ is doubly stochastic [cf. \eqref{eqn_conditions_on_weights}], which implies that $w_{ii}=1-\sum_{j\in\ccalN_i}w_{ij}<1$ as long as $w_{ij}>0$.

{The lower bound $m$ for the eigenvalues of local objective function Hessians $\nabla^2 f_i(\bbx)$ is equivalent to the strong convexity of local objective functions $f_{i}(\bbx)$ with parameter $m$. The strong convexity assumption for the local objective functions $f_{i}(\bbx)$ stated in Assumption \ref{convexity_assumption} is customary in convergence proofs of Newton-based methods, since the Hessian of objective function should be invertible to establish Newton's method [Chapter 9 of \cite{boyd}]. The upper bound $M$ for the eigenvalues of local objective function Hessians $\nabla^2 f_i(\bbx)$ is similar to the condition that gradients $\nabla f_i(\bbx)$ are Lipschitz continuous with parameter $M$ for the case that functions are twice differentiable. 

The restriction imposed by Assumption \ref{Lipschitz_assumption} is also typical of second order methods \cite{ZarghamEtal14}. Assumption \ref{Lipschitz_assumption} guarantees that the Hessian matrices of objective functions $F(\bby)$ are also Lipschitz continuous as we show in the following lemma.

%
\begin{lemma}\label{Hessian_Lipschitz_countinous}
Consider the definition of objective function $F(\bby)$ in \eqref{centralized_opt_problem}. If Assumption \ref{Lipschitz_assumption} holds then the objective function Hessian $\bbH(\bby) =: \nabla^2F(\bby)$ is Lipschitz continuous with parameter $\alpha L$, i.e.
\begin{equation}\label{H_Lipschitz_claim}
\left\|\bbH(\bby)-\bbH(\hby)\right\| \leq \alpha L \| \bby-\hby\| .
\end{equation}
for all $\bby,\hby\in \reals^{np}$.
\end{lemma}

%
\begin{myproof}
See Appendix \ref{app_Hessian_Lipschitz}.
\end{myproof}

%
Lemma \ref{Hessian_Lipschitz_countinous} states that the penalty objective function introduced in  \eqref{centralized_opt_problem} has the property that the Hessians are Lipschitz continuous, while the Lipschitz constant is a function of the penalty coefficient $1/\alpha$. This observation implies that as we increase the penalty coefficient $1/\alpha$, or, equivalently, decrease $\alpha$, the objective function $F(\bby)$ approaches a quadratic form because the curvature becomes constant. 

To prove convergence properties of NN we need bounds for the eigenvalues of the block diagonal matrix $\bbD_t$, the block sparse matrix $\bbB$, and the Hessian $\bbH_t$. These eigenvalue bounds are established in the following proposition using the conditions imposed by Assumptions \ref{ass_weight_bounds} and \ref{convexity_assumption}.

%
\begin{proposition}\label{eigenvalue_bounds}
Consider the definitions of matrices $\bbH_t$, $\bbD_t$, and $\bbB$ in \eqref{Hessian}, \eqref{diagonal_matrix}, and \eqref{non_diagona_matrix}, respectively. If Assumptions \ref{ass_weight_bounds} and \ref{convexity_assumption} hold true, then the eigenvalues of matrices $\bbH_t$, $\bbD_t$, and $\bbB$ are uniformly bounded as 
\begin{align}
 \alpha m \bbI\ \preceq \ & {\bbH_t} \ \preceq \ (2(1-\delta) + \alpha M)\bbI ,\\
 (2(1-\Delta)+\alpha m) \bbI\ \preceq\ & {\bbD_t}\ \preceq\ (2(1-\delta)+\alpha M)\bbI ,\\
 \bb0\ \preceq\ & \ {\bbB}\ \preceq\ 2 (1-\delta) \bbI.
\end{align}
\end{proposition}

%
\begin{myproof} 
See Appendix \ref{app_eigenvalue_bounds}.  
\end{myproof}

%
Proposition \ref{eigenvalue_bounds} states that Hessian matrix $\bbH_t$ and block diagonal matrix $\bbD_t$ are positive definite, while matrix $\bbB$ is positive semidefinite. 

As we noted in Section \ref{subsec_nn}, for the expansion in \eqref{exact_Hessian_inverse} to be valid the eigenvalues of the matrix $\bbD_t^{-{1}/{2}}\bbB\bbD_t^{-{1}/{2}}$ must be nonnegative and strictly smaller than $1$. The following proposition states that this is true for all times $t$.

%
\begin{proposition}\label{symmetric_term_bounds11}
Consider the definitions of the matrices $\bbD_t$ in \eqref{diagonal_matrix} and $\bbB$ in \eqref{non_diagona_matrix}. If Assumptions \ref{ass_weight_bounds} and \ref{convexity_assumption} hold true, the matrix $\bbD_t^{-{1}/{2}}  \bbB   \bbD_t^{-{1}/{2}}$ is positive semidefinite and its eigenvalues are bounded above by a constant $\rho<1$
\begin{equation}\label{important_claim}
    \bb0\  \preceq \  \bbD_t^{-{1}/{2}}  \bbB   \bbD_t^{-{1}/{2}} \ \preceq\ \rho\bbI ,
\end{equation}
where $\rho:= 2(1-\delta)/(2(1-\delta)+{\alpha m} )$.
\end{proposition}

%
\begin{myproof}
See Appendix \ref{app_eigenvalue_of_DBD}.
\end{myproof} 

%
The bounds for the eigenvalues of matrix $\bbD_t^{-{1}/{2}}  \bbB   \bbD_t^{-{1}/{2}}$ in \eqref{important_claim} guarantee convergence of the Taylor series in \eqref{exact_Hessian_inverse}. As mentioned in Section \ref{subsec_nn}, NN-$K$ truncates the first $K$ summands of the Hessian inverse Taylor series in \eqref{exact_Hessian_inverse} to approximate the Hessian inverse of the objective function in optimization problem \eqref{centralized_opt_problem}. To evaluate the performance of NN-$K$ we study the error of the Hessian inverse approximation by defining the \textit{error matrix} $\bbE_t\in\reals^{np\times np}$ as  
\begin{equation}\label{error_matrix}
\bbE_t:= \bbI- {\hbH_{t}^{(K)^{-{1}/{2}}}} \bbH_t  {\hbH_{t}^{(K)^{-{1}/{2}}}}.
\end{equation}
Error matrix $\bbE_t$ measures closeness of the Hessian inverse approximation matrix $\hbH_t^{(K)^{-1}}$ and the exact Hessian inverse $\bbH^{-1}_t$ at time $t$. Based on the definition of error matrix $\bbE_t$, if the Hessian inverse approximation $\hbH_t^{(K)^{-1}}$ approaches the exact Hessian inverse $\bbH_t^{-1}$ the error matrix $\bbE_t$ approaches the zero matrix $\bb0$. We therefore bound the error of the Hessian inverse approximation by developing a bound for the eigenvalues of the error matrix $\bbE_t$. This bound is provided in the following proposition where we further show that the error of the Hessian inverse approximation for NN-$K$ decreases exponentially as we increases $K$.

%
\begin{proposition}\label{error_matrix_proposition}
Consider the NN-$K$ method as introduced in \eqref{diagonal_matrix}-\eqref{update_formula_NN} and the definition of error matrix $\bbE_t$ in \eqref{error_matrix}. Further, recall the definition of the constant $\rho:=2(1-\delta)/(\alpha+2(1-\delta))<1$ in Proposition \ref{symmetric_term_bounds11}. The error matrix $\bbE_t$ is positive semidefinite and all its eigenvalues are upper bounded by $\rho^{K+1}$,  
\begin{equation}\label{bound_for_error}
    \bb0\  \preceq\ \bbE_t\ \preceq\ \rho^{K+1}\bbI.
\end{equation}
\end{proposition}

%
\begin{myproof} See Appendix \ref{app_error_matrix}. \end{myproof}

%
Proposition \ref{error_matrix_proposition} asserts that the error in the approximation of the Hessian inverse, thereby on the approximation of the Newton step, is bounded by $\rho^{K+1}$. This result corroborates the intuition that the larger $K$ is, the closer that $\bbd_{i,t}^{(K)}$ approximates the Newton step. This closer approximation comes at the cost of increasing the communication cost of each descent iteration. The decrease of this error being proportional to $\rho^{K+1}$ hints that using a small value of $K$ should suffice in practice. This has been corroborated in numerical experiments where $K=1$ and $K=2$ tend to work best -- see {\cite{NN-part2}}. Further note that to decrease $\rho$ we can increase $\delta$ or increase $\alpha$. Increasing $\delta$ calls for assigning substantial weight to $w_{ii}$. Increasing $\alpha$ comes at the cost of moving the solution of \eqref{centralized_opt_problem} away from the solution of \eqref{original_optimization_new_notation} and its equivalent \eqref{original_optimization_problem1}.

Bounds on the eigenvalues of the objective function Hessian $\bbH_t$ are central to the convergence analysis of Newton's method [Chapter 9 of\cite{boyd}]. Lower bounds for the Hessian eigenvalues guarantee that the matrix is nonsingular. Upper bounds imply that the minimum eigenvalue of the Hessian inverse $\bbH^{-1}$ is strictly larger than zero, which, in turn, implies a strict decrement in each Newton step. Analogous bounds for the eigenvalues of the NN approximate Hessian inverses $ {\hbH_{t}^{(K)^{-1}}}$ are required. These bounds are studied in the following lemma.

%
\begin{lemma}\label{Hessian_inverse_eigenvalue_bounds_lemma}
Consider the NN-$K$ method as defined in \eqref{diagonal_matrix}-\eqref{update_formula_NN}. If Assumptions \ref{ass_weight_bounds} and \ref{convexity_assumption} hold true, the eigenvalues of the approximate Hessian inverse $\hbH_t^{(K)^{-1}}$ are bounded as
\begin{equation}\label{bounded_Hessian_inverse}
\lambda \bbI\ \preceq\  \hbH_t^{(K)^{-1}}  \preceq\   \Lambda \bbI,
\end{equation} 
where constants $\lambda$ and $\Lambda$ are defined as
\begin{equation}\label{definition_of_lambdas}
\!\!\! \lambda\!:=\! \frac{1}{2(1-\delta)+\alpha M } \ \text{ and} 
\ \ \Lambda\! :=\! {\frac{1-\rho^{K+1}}{(1-\rho)(2(1-\Delta)+\alpha m )}}.
\end{equation}
\end{lemma}

%
\begin{myproof}
See Appendix \ref{app_eigenvalues_of_Hessian__inverse_approximation}.
\end{myproof}

%
According to the result of Lemma \ref{Hessian_inverse_eigenvalue_bounds_lemma}, the NN-$K$ approximate Hessian inverses $\hbH_t^{(K)^{-1}}$ are strictly positive definite and have all of their eigenvalues bounded between the positive and finite constants $\lambda$ and $\Lambda$. This is true for all $K$ and uniform across all iteration indexes $t$. Considering these eigenvalue bounds and the fact that $-\bbg_t$ is a descent direction, the approximate Newton step $-\hbH_t^{(K)^{-1}}\bbg_t$ enforces convergence of the iterate $\bby_t$ to the optimal argument $\bby^*$ of the penalized objective function $F(\bby)$ in \eqref{centralized_opt_problem}. In the following theorem we show that if the stepsize $\epsilon$ is properly chosen, the sequence of objective function values $F(\bby_t)$ converges at least linearly to the optimal objective function value $F(\bby^*)$.

%
\begin{thm}\label{linear_convergence}
Consider the NN-$K$ method as defined in \eqref{diagonal_matrix}-\eqref{update_formula_NN} and the objective function $F(\bby)$ as introduced in \eqref{centralized_opt_problem}. Further, recall the definitions of the lower and upper bounds $\lambda$ and $\Lambda$, respectively, for the eigenvalues of the approximate Hessian inverse $ \hbH_t^{(K)^{-1}} $ in \eqref{definition_of_lambdas}. If the stepsize $\epsilon$ is chosen as 
\begin{equation}\label{step_size_condition}
\epsilon = \min  \left\{ 1\ , \left[{\frac{3m\lambda^{\frac{5}{2}}}{ L\Lambda^{3}{(F(\bby_{0})-F(\bby^{*})) }^{\frac{1}{2}} } }  \right]^{\frac{1}{2}} \right\} 
\end{equation}
and Assumptions \ref{ass_weight_bounds}, \ref{convexity_assumption}, and \ref{Lipschitz_assumption} hold true, the sequence $F(\bby_{t})$ converges to the optimal argument $F(\bby^{*})$ at least linearly with constant $0<1-\zeta<1$. I.e.,
\begin{equation}\label{linear_convegrence_claim}
F(\bby_{t}) -F(\bby^*) \leq (1-\zeta)^t  {\left(F(\bby_{0}) -F(\bby^*) \right)},
\end{equation}
where the constant $0<\zeta<1$ is explicitly given by
\begin{equation}\label{beta_0_defintion}
\zeta :=    {(2-\eps) \epsilon\alpha m\lambda} - \frac{\alpha\epsilon^3  L\Lambda^{3}(F(\bby_{0})-F(\bby^{*}) )^{\frac{1}{2}}}{6\lambda^{\frac{3}{2}}} .
\end{equation}\end{thm}

%
\begin{myproof} See Appendix \ref{app_linear_convergecene_theorem}. \end{myproof}

%
Theorem \ref{linear_convergence} shows that the objective function error sequence $F(\bby_t)-F(\bby^{*})$ asymptoticly converges to zero and that the rate of convergence is at least linear. Note that according to the definition of the convergence parameter $\zeta $ in Theorem \ref{linear_convergence} and the definitions of $\lambda$ and $\Lambda$ in \eqref{definition_of_lambdas}, increasing $\alpha$ leads to faster convergence. This observation verifies existence of a tradeoff between rate and accuracy of convergence. For large values of $\alpha $ the sequence generated by Network Newton converges faster to the optimal solution of \eqref{centralized_opt_problem}. These faster convergence comes at the cost of increasing the distance between the optimal solutions of \eqref{centralized_opt_problem} and \eqref{original_optimization_problem1}. Conversely, smaller $\alpha $ implies smaller gap between the optimal solutions of \eqref{centralized_opt_problem} and \eqref{original_optimization_problem1}, but the convergence rate of NN-$K$ is slower. This suggests value in the use of adaptive strategies for the selection of $\alpha$ that we develop in {\cite{NN-part2}.}

\section{Conclusions}

This paper developed the network Newton method as an approximate Newton method for solving distributed optimization problems where the components of the objective function are available at different nodes of a network. The algorithm builds on a reinterpretation of distributed gradient descent as a penalty method and relies on an approximation of the Newton step of the corresponding penalized objective function. To approximate the Newton direction we truncate the Taylor series of the exact Newton step. This leads to a family of methods defined by the number $K$ of Taylor series terms kept in the approximation. When we keep $K$ terms of the Taylor series, the method is called NN-$K$ and can be implemented through the aggregation of information in $K$-hop neighborhoods. We showed that the proposed method converges at least linearly to the solution of the penalized objective, and, consequently, to a neighborhood of the optimal argument for the original optimization problem. It follows from this convergence analysis that larger penalty coefficients result in faster convergence that comes at the cost of increasing the distance between the optimal solutions of the original and penalized objectives. 

This paper does not show any advantage of NN relative to distributed gradient descent, other than the expectation to see improved convergence times due to the attempt to approximate the Newton direction of the penalized objective. These advantages are shown in a companion paper where we: (i) Show that the convergence rate is quadratic in a specific interval that can be made arbitrarily large by increasing $K$ \cite{NN-part2}. (ii) Use numerical results to establish the advantages of NN-$K$ in terms of number of iterations and communications steps relative to DGD.

\begin{appendices}


\section{Proof of Lemma  \ref{Hessian_Lipschitz_countinous}}\label{app_Hessian_Lipschitz}

Consider two vectors $\bby:=[\bbx_1;\dots;\bbx_n]  \in \reals^{np}$ and $\hby:=[\hbx_1;\dots;\hbx_n] \in \reals^{np}$. Based on the Hessian expression in \eqref{Hessian}, we simplify the Euclidean norm of the Hessian difference $ \bbH(\bby)-\bbH(\hby)$ as
\begin{equation}\label{H_Lipschitz_continuity}
\|\bbH(\bby)-\bbH(\hby)\| 
				=  \alpha \left\|  \bbG(\bby)-\bbG(\hby) \right\|.
\end{equation}
The result in \eqref{H_Lipschitz_continuity} is implied by the fact that the matrix $\bbI-\bbZ$ does not depend on the argument $\bby$ of the Hessian $\bbH(\bby)$.
The next step is to bound the norm of the difference for two $\bbG$ matrices $\| \bbG(\bby)-\bbG(\hby)\|$ in terms of the difference between two vectors $\bby$ and $\hby$, i.e. $\|\bby-\hby\|$.

According to the definition of $\bbG(\bby)$ in \eqref{G_form}, the difference matrix $\bbG(\bby)-\bbG(\hby)$ is block diagonal and the $i$th diagonal block is
\begin{equation}\label{matrix_difference}
{\bbG(\bby)_{ii}-\bbG(\hby)}_{ii}=
\nabla^2 f_{i}(\bbx_{i})-\nabla^2 f_{i}(\hbx_{i}) .
\end{equation}
Consider any vector $\bbv\in \reals^{np}$ and separate each $p$ components of vector $\bbv$ and consider it as a new vector called $\bbv_{i}\in \reals^p$, i.e. $\bbv:=[\bbv_1;\dots;\bbv_n]$. Observing the relation for the difference $\bbG(\bby)-\bbG(\hby)$ in \eqref{matrix_difference}, the symmetry of matrices $\bbG(\bby)$ and $\bbG(\hby)$, and the definition $\|\bbA\|_2:=\sqrt{\lambda_{max}(\bbA^T\bbA)}$ of the Euclidean norm of a matrix, we obtain that the squared difference norm $\|\bbG(\bby)-\bbG(\hby)\|_2^2$ can be written as 
\begin{align}\label{inner_product_differnece}
\left\|\bbG(\bby)\!-\!\bbG(\hby)\right\|_{2}^2 \!
	&= \max_{\bbv} \frac{\bbv^T[\bbG(\bby)-\bbG(\hby)]^2\bbv}{\|\bbv\|^2} \nonumber \\ 
    &= \max_{\bbv} \frac{\sum_{i=1}^n \bbv_{i}^T \left[\nabla^2 f_{i}(\bbx_{i})-\nabla^2 f_{i}(\hbx_{i})\right]^2 	 \bbv_{i}}		                     {\|\bbv\|^2}
\end{align}
Observe that each summand in \eqref{inner_product_differnece} can be upper bounded by applying Cauchy-Schwarz inequality as 
\begin{equation}\label{couchy_result}
\bbv_{i}^T\! \left[\nabla^2 f_{i}(\bbx_{i})\!-\!\nabla^2 f_{i}(\hbx_{i})\right]^2\! \bbv_{i} \!\leq\!
		\left\|\nabla^2 f_{i}(\bbx_{i})\!-\!\nabla^2 f_{i}(\hbx_{i})\right\|_{2}^2 \! \|\bbv_{i}\|^2
\end{equation}
Substituting the upper bound in \eqref{couchy_result} into \eqref{inner_product_differnece} implies that the squared norm $\left\|\bbG(\bby)-\bbG(\hby)\right\|_{2}^2$ is bounded above as
\begin{equation}\label{norm_difference_2}
\!\left\|\bbG(\bby)\!-\!\bbG(\hby)\right\|_{2}^2\leq
 \max_{\bbv} \frac{\sum_{i=1}^n\!\left\|\nabla^2 f_{i}(\bbx_{i})\!-\!\nabla^2 f_{i}(\hbx_{i})\right\|_{2}^2  \|\bbv_{i}\|^2}{\|\bbv\|^2}.
\end{equation}
Observe that Assumption 2 states that the local objective function Hessians $\nabla^2 f_{i}(\bbx_{i})$ are Lipschitz continuous with parameter $L$, i.e., $\|\nabla^2 f_{i}(\bbx_{i})-\nabla^2 f_{i}(\hbx_{i})\|\leq L \|\bbx_{i}-\hbx_{i}\|$. Considering this inequality the upper bound in \eqref{norm_difference_2} can be changed by replacing $\|\nabla^2 f_{i}(\bbx_{i})-\nabla^2 f_{i}(\hbx_{i})\|$ by $L \|\bbx_{i}-\hbx_{i}\|$ which yields
\begin{equation}\label{alef}
\left\|\bbG(\bby)-\bbG(\hby)\right\|_{2}^2\leq
 \max_{\bbv} \frac{L^2\sum_{i=1}^n\left\|\bbx_{i}-\hbx_{i}\right\|_{2}^2  \|\bbv_{i}\|^2}{\sum_{i=1}^n\|\bbv_{i}\|^2}.
\end{equation}
Note now that for any sequences of scalars $a_{i}$ and $b_i$, the inequality $\sum_{i=1}^n a_{i}^2b_i^2\leq (\sum_{i=1}^n a_{i}^2)(\sum_{i=1}^n b_{i}^2)$ holds. If we divide both sides of this relation by $\sum_{i=1}^n b_{i}^2$ and set $a_i=\|\bbx_{i}-\hbx_{i}\|$ and $b_i=\|\bbv_{i}\|$, we obtain
\begin{equation}\label{alef_2}
\frac{\sum_{i=1}^n\left\|\bbx_{i}-\hbx_{i}\right\|_{2}^2  \|\bbv_{i}\|^2}{\sum_{i=1}^n\|\bbv_{i}\|^2}
	\leq \sum_{i=1}^n\left\|\bbx_{i}-\hbx_{i}\right\|_{2}^2.
\end{equation}
Combining the two inequalities in \eqref{alef} and \eqref{alef_2} leads to
\begin{equation}\label{alef_3}
\left\|\bbG(\bby)-\bbG(\hby)\right\|_{2}^2
   \leq \max_{\bbv} L^2 \sum_{i=1}^n\left\|\bbx_{i}-\hbx_{i}\right\|_{2}^2  .
\end{equation}
Since the right hand side of \eqref{alef_3} does not depend on the vector $\bbv$ we can eliminate maximization with respect $\bbv$. Further, note that according to the structure of vectors $\bby$ and $\hby$, we can write $\left\|\bby-\hby\right\|_{2}^2=\sum_{i=1}^n\left\|\bbx_{i}-\hbx_{i}\right\|_{2}^2$. These two observations in association with \eqref{alef_3} imply that squared norm of the difference between matrices $\bbG(\bby)$ and $\bbG(\hby)$ is bounded above by
\begin{equation}\label{alef_2223}
\left\|\bbG(\bby)-\bbG(\hby)\right\|_{2}^2 \leq  L^2 \left\|\bby-\hby\right\|_{2}^2.
\end{equation}
Taking the square root of both sides of \eqref{alef_2223} implies 
\begin{equation}\label{Lipschitz_countinuity_of_G}
\left\|\bbG(\bby)-\bbG(\hby)\right\|_{2} \leq  {L} \left\|\bby-\hby\right\|_{2}.
\end{equation}
According to \eqref{Lipschitz_countinuity_of_G} we can conclude that the matrix $\bbG$ is Lipschitz continuos with parameter $L$. Considering the expression in \eqref{H_Lipschitz_continuity} and the inequality in \eqref{Lipschitz_countinuity_of_G}, the claim in \eqref{H_Lipschitz_claim} follows.



\section{Proof of Proposition  \ref{eigenvalue_bounds}}\label{app_eigenvalue_bounds}

We first study the bounds for the eigenvalues of matrix $\bbI-\bbZ$. Notice that since $\bbI-\bbZ$ can be written as $( \bbI_n - \bbW) \otimes \bbI_p$, all the eigenvalues of matrix $\bbI-\bbZ$ are in the spectrum of matrix $\bbI-\bbW$. Therefore, we can study the bounds for the eigenvalues of matrix $\bbI-\bbW$ in lieu of matrix $\bbI-\bbZ$. The Gershgorin circle theorem states that each eigenvalue of a matrix $\bbA$ lies within at least one of the Gershgorin discs $D(a_{ii}, R_{ii})$ where the center $a_{ii}$ is the $i$th diagonal element of $A$ and the radius $R_{ii}:=\sum_{j \neq i} |a_{ij}|$ is the sum of the absolute values of all the non-diagonal elements of the $i$th row. Note that matrix $\bbI-\bbW$ is symmetric and as a result all eigenvalues are real. Hence, Gershgorin discs can be considered as intervals of width $[a_{ii}-R_{ii},a_{ii}+R_{ii}]$ for matrix $\bbI-\bbW$, where $a_{ii}=1-w_{ii}$ and $R_{ii}= \sum_{j \neq i} |w_{ij}|$. Since all the elements of matrix $\bbW$ are non-negative, $|w_{ij}|$ can be substituted by $w_{ij}$. Therefore, all the eigenvalues of matrix $\bbI-\bbW$ in at least one of the intervals $[1-w_{i} - \sum_{j \neq i} w_{ij},1-w_{i} + \sum_{j \neq i} w_{ij}]$. Now observing that sum of the weights that a node assigns to itself and all the other nodes is one, i.e. $\sum_{j} w_{ij}=1$, it can be derived that $1-w_{ii} =\sum_{j \neq i}^{n} w_{ij} $. This observation implies that the Gershgorin intervals can be simplified as $[0,2(1-w_{ii})]$ for $i=1,\dots, n$. This observation in association with the fact that $2(1-w_{ii})\leq 2(1-\delta)$ implies that all the eigenvalues of matrix $\bbI-\bbW$ are in the interval $[0,2(1-\delta)]$ and consequently the eigenvalues of matrix $\bbI-\bbZ$ are bounded as
\begin{equation}\label{I-Z_bounds}
\bb0\ \preceq\ \bbI-\bbZ\ \preceq \ 2(1-\delta) \bbI.
\end{equation} 
To prove bounds for the eigenvalues of Hessian $\bbH_t$, first we find lower and upper bounds for the eigenvalues of matrix $\bbG_t$. Since matrix $\bbG_t$ is block diagonal and the eigenvalues of each diagonal block $\bbG_{ii,t}=\nabla^2 f_{i}(\bbx_{i,t})$ are bounded by constants $0<m\leq M<\infty$ as mentioned in \eqref{local_hessian_eigenvlaue_bounds}, we obtain that the eigenvalues of matrix $\bbG_t$ are bounded as 
\begin{equation}\label{G_bounds}
m\bbI\ \preceq\ \bbG_t\ \preceq \ M\bbI.
\end{equation} 
Considering the definition of Hessian $\bbH_t:=\bbI-\bbZ+\alpha \bbG_t$ and the bounds in \eqref{I-Z_bounds} and \eqref{G_bounds}, the first claim follows.

We proceed now to prove bounds for the eigenvalues of block diagonal matrix $\bbD_t$. According to the definition of matrix $\bbD_t$ in \eqref{diagonal_matrix} we can write 
\begin{equation}\label{D_new_formulation}
\bbD_{t}\ =\ \alpha \bbG_t \ + \ (\bbI_n - \bbW_{d})\otimes \bbI_p\ ,
\end{equation}
where $\bbW_d$ is defined as $\bbW_{d}:=\diag(\bbW)$. Note that matrix $\bbI_n - \bbW_{d}$ is diagonal and the $i$-th diagonal component is $1-w_{ii}$. Since the local weights satisfy $\delta \leq w_{ii} \leq \Delta$, we obtain that eigenvalues of matrix $\bbI_n - \bbW_{d}$ are bounded below and above by $1-\Delta $ and $1-\delta$, respectively. Observe that the eigenvalue sets of matrices $(\bbI_n - \bbW_{d})$ and $(\bbI_n - \bbW_{d})\otimes \bbI_p$ are identical which implies 
\begin{equation}\label{D_second_part_bounds}
 (1-\Delta ) \bbI_{np} \   \preceq  \ (\bbI_n - \bbW_{d})\otimes \bbI_p  \ \preceq \
 	(1-\delta ) \bbI_{np}
\end{equation}
Considering the relation in \eqref{D_new_formulation} and bounds in \eqref{G_bounds} and \eqref{D_second_part_bounds}, the second claim follows.

Based on the definition of matrix $\bbB$ in \eqref{non_diagona_matrix} and the relation that $\bbZ=\bbW\otimes \bbI$ we can write
\begin{equation}\label{simplificaion_for_B}
\bbB\ 
=\ (\bbI - 2 \bbW_{d} + \bbW)  \otimes  \bbI .
\end{equation}
Hence, to bound eigenvalues of matrix $\bbB$ we study lower and upper bounds for the eigenvalues of matrix $\bbI - 2 \bbW_{d} + \bbW$. Observe that in the $i$-{th} row of matrix $\bbI - 2 \bbW_{d} + \bbW$, the  diagonal component is $1-w_{ii}$ and the $j$th component is $ w_{ij} $ for all $j\neq i$. Using Gershgorin theorem and the same argument that we established for the eigenvalues of $\bbI - \bbZ$, we can write 
\begin{equation}\label{last_step_of_proof}
 \bb0\ \preceq\ \bbI - 2 \bbW_{d} + \bbW\  \preceq\ 2(1-\delta)\bbI.
\end{equation}
Considering \eqref{last_step_of_proof} and the expression for matrix $\bbB$ in \eqref{simplificaion_for_B}, the last claim follows.


\section{Proof of Proposition  \ref{symmetric_term_bounds11}}\label{app_eigenvalue_of_DBD}

According to the result of Proposition 1, the block diagonal matrix $\bbD_{t}$ is positive definite and matrix $\bbB$ is positive semidefinite which immediately implies that matrix $\bbD_t^{-{1}/{2}}\bbB\bbD_t^{-{1}/{2}}$ is positive semidefinite and the lower bound in \eqref{important_claim} follows. 

Recall the definition of block diagonal matrix $\bbD_t$ in \eqref{diagonal_matrix} and define matrix $\hbD$ as a special case of matrix $\bbD_t$ for $\alpha=0$. I.e., $\hbD:=2(\bbI-\bbZ_d)$. Notice that matrix $\hbD$ is diagonal, only depends on the structure of the network, and that it is also time invariant. Since matrix $\hbD$ is diagonal and each diagonal component $1-w_{ii}$ is strictly larger than 0, matrix $\hbD$ is positive definite and invertible. Therefore, we can write $\bbD_t^{-{1}/{2}}\bbB\bbD_t^{-{1}/{2}}$ as
\begin{equation}\label{new_decomposition}
\bbD_t^{-\frac{1}{2}}\bbB\bbD_t^{-\frac{1}{2}}=
	\left(	\bbD_t^{-\frac{1}{2}} \hbD^{\frac{1}{2}}	\right)
	\left(	\hbD^{-\frac{1}{2}}\bbB\hbD^{-\frac{1}{2}}	\right)
	\left(	\hbD^{\frac{1}{2}}\bbD_t^{-\frac{1}{2}}	\right).
\end{equation}
The next step is to find an upper bound for the eigenvalues of the symmetric term $\hbD^{-{1/2}}\bbB\hbD^{-{1/2}}$ in \eqref{new_decomposition}. Observing the fact that matrices $\hbD^{-{1}/{2}}\bbB\hbD^{-{1}/{2}}$ and $\bbB\hbD^{-1}$ are \textit{similar}, eigenvalues of these matrices are identical. Therefore, we proceed to characterize an upper bound for the eigenvalues of matrix $\bbB\hbD^{-1}$. Based on the definitions of matrices $\bbB$ and $\hbD$, the product $\bbB\hbD^{-1}$ is given by
\begin{equation}\label{product_of_B_and_D_hat}
\bbB\hbD^{-1} = \left(\bbI-2\bbZ_{d}+\bbZ \right) (2(\bbI-\bbZ_{d}))^{-1}.
\end{equation}
Therefore, the general form of matrix $\bbB\hbD^{-1}$ is
\begin{equation}\label{B_D_hat_inverse_form}
\bbB\hbD^{-1}= 
	\frac{1}{2}\left[ {\begin{array}{*{20}c}
\bbI &
\frac{w_{12}}{(1-w_{22})}\bbI& \ldots &
\frac{w_{1n}}{(1-w_{nn})}\bbI \\ 
	\frac{w_{21}}{(1-w_{11})}\bbI &\bbI  &\ldots & 
	\frac{w_{2n}}{(1-w_{nn})}\bbI\\
	 \vdots & \vdots & \ddots & \vdots \\ 
	\frac{w_{n1}}{(1-w_{11})}\bbI &\frac{w_{n2}}{(1-w_{22})}\bbI& \ldots &
	\bbI
 \end{array} } \right].
\end{equation}
Note that each diagonal component of matrix $\bbB\hbD^{-1}$ is $1/2$ and that the sum of non-diagonal components of column $i$ is 
\begin{equation}\label{sum_of_non_diagonal_components}
   \sum_{j=1,j\neq i}^{np} \bbB\hbD^{-1}[ji]
       = \frac{1}{2} \sum_{j=1,j\neq i}^{np} \frac{w_{ji}}{1-w_{ii}}
       = \frac{1}{2}.
\end{equation}
Now, by considering the result in \eqref{sum_of_non_diagonal_components} and applying Gershgorin theorem we can conclude that eigenvalues of matrix $\bbB\hbD^{-1}$ are bounded as
\begin{equation}\label{bounds_for_eigenvalues_of_B_D_hat_^-1}
0 \leq \mu_{i}(\bbB\hbD^{-1}) \leq 1      \qquad i=1,\dots, n,
\end{equation}
where $\mu_{i}(\bbB\hbD^{-1})$ indicates the $i$-{th} eigenvalue of matrix $\bbB\hbD^{-1}$. The bounds in \eqref{bounds_for_eigenvalues_of_B_D_hat_^-1}  and \textit{similarity} of matrices $\bbB\hbD^{-1}$ and $\hbD^{-1/2}\bbB\hbD^{-1/2} $ show that the eigenvalues of matrix $\hbD^{-1/2}\bbB\hbD^{-1/2} $ are uniformly bounded in the interval 
\begin{equation}\label{bounds_for_eigenvalues_of_D_B_D}
0 \leq \mu_{i}(\hbD^{-{1}/{2}}\bbB\hbD^{-{1}/{2}}) \leq 1.
\end{equation}

Based on the decomposition in \eqref{new_decomposition} to characterize the bounds for the eigenvalues of matrix $\bbD_t^{-1/2}\bbB\bbD_t^{-1/2}$, the bounds for the eigenvalues of matrix $\hbD^{1/2}\bbD_t^{-1/2}$ should be studied as well. Notice that according to the definitions of matrices $\hbD$ and $\bbD_t$, the product $\hbD^{1/2}\bbD_t^{-1/2}$ is block diagonal and the $i$-th diagonal block is 
\begin{equation}\label{product_DDD}
\left[{\hbD^{{1}/{2}} \bbD_{t}^{-{1}/{2}}}\right]_{ii}= 
\left(\frac{ \alpha  \nabla^2 f_{i}(\bbx_{i,t})}{2(1-w_{ii})} +\bbI\right)^{-1/2}.
\end{equation}
Observe that according to Assumption 1, the eigenvalues of local Hessian matrices $\nabla^{2}f_{i}(\bbx_{i})$ are bounded by $m$ and $M$. Further notice that the diagonal elements of weight matrix $w_{ii}$ are bounded by $\delta $ and $\Delta $, i.e. $\delta\leq w_{ii} \leq \Delta$. Considering these bounds we can show that the eigenvalues of matrices $ (\alpha/2(1-w_{11}))  \nabla^2 f_{i}(\bbx_{i,t})+\bbI $ are lower and upper bounded as 
\begin{equation}\label{bounds_33}
\left[\frac{\alpha m}{2(1-\delta)}+1\right]  \bbI  
\preceq  \frac{ \alpha  \nabla^2 f_{i}(\bbx_{i,t})}{2(1-w_{11})}+\bbI 
\preceq \left[\frac{\alpha M}{2(1-\Delta)}+1\right]  \bbI .
\end{equation}
By considering the bounds in \eqref{bounds_33}, the eigenvalues of each block of matrix $\hbD^{1/2}\bbD_t^{-1/2}$ as introduced in \eqref{product_DDD} are bounded below and above as
\begin{align}\label{bounds_3}
\left[\frac{2(1-\Delta)}{2(1-\Delta)+\alpha M} \right]^{\frac{1}{2}}  \bbI  
&\preceq   \left[\frac{ \alpha  \nabla^2 f_{i}(\bbx_{i,t})}{2(1-w_{11})}+\bbI \right]^{-\frac{1}{2}} \\
&\qquad \qquad \preceq  \left[ \frac{2(1-\delta)}{2(1-\delta)+\alpha m}  \right]^{\frac{1}{2}}  \bbI.
\nonumber 
\end{align}
Since \eqref{bounds_3} holds for all the diagonal blocks of matrix $\hbD^{1/2}\bbD_t^{-1/2}$, the eigenvalues of this matrix also satisfy the bounds in \eqref{bounds_3} which implies that
\begin{equation}\label{bounds_for_DDD}
	\left[\frac{2(1-\Delta)}{2(1-\Delta)+\alpha M} \right]^{\frac{1}{2}}\!\!
		\leq \mu_{i}\!\left( \hbD^{\frac{1}{2}} \bbD_t^{-\frac{1}{2}} \right) 
		\leq 	\left[ \frac{2(1-\delta)}{2(1-\delta)+\alpha m}\right]^{\frac{1}{2}}\!\!,
\end{equation}
for $i=1, \dots, n$.
Observing the decomposition in \eqref{new_decomposition}, the norm of the matrix $\bbD_t^{-{1}/{2}}\bbB\bbD_t^{-{1}/{2}}$ is upper bounded as
\begin{equation}\label{new_decomposition_2}
\!\!\left\|\bbD_t^{-{1}/{2}}\bbB\bbD_t^{-{1}/{2}} \right\|\leq 
	\left\|	\bbD_t^{-{1}/{2}} \hbD^{{1}/{2}}	\right\|^2
	\left\|	\hbD^{-{1}/{2}}\bbB\hbD^{-{1}/{2}}	\right\|.
\end{equation}
Considering the symmetry of matrices $\hbD^{1/2}\bbD_t^{-1/2}$ and $\hbD^{-{1}/{2}}\bbB\hbD^{-{1}/{2}}$, and the upper bounds for their eigenvalues in \eqref{bounds_for_eigenvalues_of_D_B_D} and \eqref{bounds_for_DDD}, respectively, we can substitute the norm of these two matrices by the upper bounds of their eigenvalues and simplify the upper bound in \eqref{new_decomposition_2} to
\begin{equation}\label{new_decomposition_3}
\left\|\bbD_t^{-{1}/{2}}\bbB\bbD_t^{-{1}/{2}} \right\|\leq \frac{2(1-\delta)}{2(1-\delta)+\alpha m}.
\end{equation}
Based on the upper bound for the norm of the matrix $\bbD_t^{-{1}/{2}}\bbB\bbD_t^{-{1}/{2}}$ in \eqref{new_decomposition_3} and the fact that matrix $\bbD_t^{-{1}/{2}}\bbB\bbD_t^{-{1}/{2}}$ is positive semidefinite, we can conclude that the eigenvalues of matrix $\bbD_t^{-{1}/{2}}\bbB\bbD_t^{-{1}/{2}}$ are upper bounded by ${2(1-\delta)}/({2(1-\delta)+\alpha m})$ and the right hand side of \eqref{important_claim} follows.


\section{Proof of Proposition  \ref{error_matrix_proposition}}\label{app_error_matrix}

In this proof and the rest of the proofs we denote the Hessian approximation as $\hbH_{t}^{-1}$ instead of $\hbH_{t}^{{(K)^{-1}}}$ for simplification of equations. To prove lower and upper bounds for the eigenvalues of the error matrix $\bbE_{t}$ we first develop a simplification for the matrix $  \bbI -\bbH_t\hbH_{t}^{-1} $ in the following lemma.

\begin{lemma}\label{error_matrix_simplfication}
Consider the NN-$K$ method as defined in \eqref{diagonal_matrix}-\eqref{update_formula_NN}. The matrix $\bbI -\bbH_t\hbH_{t}^{-1} $ can be simplified as
\begin{equation}\label{error_lemma_claim}
\bbI -\bbH_t\hbH_{t}^{-1} =  \left(   \bbB    \bbD_{t}^{-1}     \right)^{K+1} .
\end{equation} 

\end{lemma}

\begin{myproof}
Considering the definitions of the Hessian inverse approximation $\hbH_{t}^{-1}$ in \eqref{Hessian_inverse_approximation} and the matrix decomposition for the exact Hessian $\bbH_t=\bbD_t-\bbB$, we obtain 
\begin{align}\label{error_1}
  \bbH_t\hbH_{t}^{-1} 
&=  \left( \bbD_t-\bbB \right)
\left(  \bbD_t^{-\frac{1}{2}}  \  \sum_{k=0}^{K} \left(  \bbD_t^{-\frac{1}{2}}  \bbB   \bbD_t^{-\frac{1}{2}}            \right)^{k}     \ \bbD_t^{-\frac{1}{2}}  \right).\nonumber\\
&=  \left( \bbD_t-\bbB \right)  \sum_{k=0}^{K} \bbD_t^{-1}\left(    \bbB \bbD_t^{-1}       \right)^{k} .
\end{align}
By considering the result in \eqref{error_1}, we simplify the expression $  \bbI -\bbH_t\hbH_{t}^{-1} $ as
\begin{align}\label{error_2}
\bbI -\bbH_t\hbH_{t}^{-1} & =  \bbI - (\bbD_t-\bbB)  \sum_{k=0}^{K}  \bbD_t^{-1} \left(  \bbB\bbD_t^{-1}        \right)^{k}    \nonumber \\
& =   \bbI -   \sum_{k=0}^{K} \left(   \bbB    \bbD_t^{-1}     \right)^{k}
     +  \sum_{k=0}^{K} \left(   \bbB    \bbD_t^{-1}     \right)^{k+1}  .
\end{align}
In the right hand side of \eqref{error_2} the identity matrix cancels out the first term in the sum $\sum_{k=0}^{K} \left(\bbB\bbD_t^{-1}\right)^{k}$. The remaining terms of this sum are cancelled out by the first $K$ terms of the sum $\sum_{k=0}^{K} \left(   \bbB    \bbD_t^{-1}     \right)^{k+1}$ so that the whole expression simplifies to $(\bbB\bbD_{t}^{-1} )^{K+1}$ as is claimed in \eqref{error_lemma_claim}.
\end{myproof}

Observing the fact that the error matrix $\bbE_t$ is a conjugate of the matrix $\bbI -\bbH_t\hbH_{t}^{-1} $ and considering the simplification in Lemma \ref{error_matrix_simplfication} we show that the eigenvalues of error matrix $\bbE_t$ are bounded.\\

\textbf{Proof of Proposition \ref{error_matrix_proposition}:}
Recall the result of Proposition \ref{symmetric_term_bounds11} that all the eigenvalues of matrix $\bbD_t^{-1/2}\bbB\bbD_t^{-1/2}$ are uniformly bounded between $0$ and $\rho$. Since matrices $\bbD_t^{-1/2}\bbB\bbD_t^{-1/2}$ and $\bbB_t    \bbD_{t}^{-1}$ are similar (conjugate) the sets of eigenvalues of these two matrices are identical. Therefore, eigenvalues of matrix $\bbB\bbD^{-1}$ are bounded as
\begin{equation}\label{bound_help}
  0 \ \leq\ \mu_{i}(\bbB\bbD^{-1})\ \leq\  \rho,
\end{equation}
for $i = {1,2,\dots,np}$. The bounds for the eigenvalues of matrix $\bbB\bbD^{-1}$ in association with expression \eqref{error_lemma_claim} leads to the following bounds for the eigenvalues of matrix $ \bbI -\bbH_t\hbH_{t}^{-1}$,
\begin{equation}\label{bound_help_2}
  0 \ \leq\ \mu_{i}(\bbI -\bbH_t\hbH_{t}^{-1})\ \leq\  \rho^{K+1}.
\end{equation}
Observe that the error matrix $\bbE_t=\bbI-\hbH_{t}^{-1/2}\bbH_t\hbH_{t}^{-1/2}$ is the conjugate of matrix $\bbI -\bbH_t\hbH_{t}^{-1}$. Hence, the bounds for the eigenvalues of matrix $\bbI -\bbH_t\hbH_{t}^{-1}$ also hold for the eigenvalues of error matrix $\bbE_t$ and the claim in \eqref{bound_for_error} follows.


\section{Proof of Lemma  \ref{Hessian_inverse_eigenvalue_bounds_lemma}}\label{app_eigenvalues_of_Hessian__inverse_approximation}

According to the Cauchy-Schwarz inequality, the product of the norms is larger than norm of the products. This observation and the definition of the approximate Hessian inverse $\hbH_{t}^{-1}$ in \eqref{Hessian_inverse_approximation} leads to
\begin{align}\label{upper_bound_for_norm_of_Hessian_inverse}
\Big{\|}\hbH_{t}^{-1}\Big{\|} \leq  \left\| \bbD_t^{-{1}/{2}}\right\|^2  &\Big\|\ \bbI + \bbD_t^{-{1}/{2}} \bbB\bbD_t^{-{1}/{2}}  + \ldots \nonumber\\
& \qquad+ \left[\bbD_t^{-{1}/{2}} \bbB\bbD_t^{-{1}/{2}}  \right]^K   \Big\|.
\end{align}
Observe that as a result of Proposition \ref{eigenvalue_bounds} the eigenvalues of matrix $\bbD_t$ are bounded below by $ {2(1-\Delta)+\alpha m }$. Therefore, the maximum eigenvalue of its inverse $\bbD_t^{-1}$ is smaller than $ 1/(2(1-\Delta)+\alpha m)$. It then follows that the norm of the matrix $\bbD_t^{-1/2}$ is bounded above as
\begin{equation}\label{up1}
\left\| \bbD_t^{-1/2} \right\|\ \leq \ \left[\frac{1}{2(1-\Delta)+\alpha m}\right]^{1/2}.
\end{equation}
Based on the result in Proposition \ref{symmetric_term_bounds11} the eigenvalues of the matrix $\bbD_t^{-{1}/{2}} \bbB\bbD_t^{-{1}/{2}}$ are smaller than $\rho$. Further using the symmetry and positive definiteness of the matrix $\bbD_t^{-{1}/{2}} \bbB\bbD_t^{-{1}/{2}}$, we obtain 
\begin{equation}\label{up2}
\left\|\bbD_t^{-{1}/{2}} \bbB\bbD_t^{-{1}/{2}} \right\| \ \leq  \ \rho.
\end{equation}
Using the triangle inequality in \eqref{upper_bound_for_norm_of_Hessian_inverse} to claim that the norm of the sum is smaller than the sum of the norms and substituting the upper bounds in \eqref{up1} and \eqref{up2} in the resulting expression we obtain 
\begin{equation}\label{up3}
\left\| \hbH_{t}^{-1} \right\|    \leq \frac{1}{2(1-\Delta)+\alpha m }\
		\sum_{k=0}^K \rho^k.
\end{equation}
By considering the fact that $\rho$ is smaller than $1$, the sum $\sum_{k=0}^K \rho^k$ can be simplified to $(1-\rho^{K+1})/(1-\rho)$. Considering this simplification for the sum in \eqref{up3}, the upper bound in \eqref{bounded_Hessian_inverse} for the eigenvalues of the approximate Hessian inverse $\hbH_{t}^{-1}$ follows.

The next step is to provide a lower bound for the eigenvalues of the Hessian inverse approximation matrix $\hbH_{t}^{-1}$. In the Hessian inverse approximation formula \eqref{Hessian_inverse_approximation}, all the summands except the first one, $\bbD_t^{-1}$, are positive semidefinite. Hence, the approximate Hessian inverse $\hbH_{t}^{-1}$ is the sum of matrix $\bbD_t^{-1}$ and $K$ positive semidefinite matrices and as a result we can conclude that 
\begin{equation}\label{lower_bounded_matrix}
\bbD_t^{-1} \preceq\  \hbH_{t}^{-1} .
\end{equation}
Proposition \ref{eigenvalue_bounds} shows that the eigenvalues of matrix $\bbD_t$ are bounded above by $2(1-\delta)+\alpha M$ which leads to the conclusion that there exits a lower bound for the eigenvalues of matrix $\bbD_t^{-1}$,
\begin{equation}\label{lower_bounded_matrix2}
\frac{1}{2(1-\delta)+\alpha M}\ \bbI \ \preceq\ \bbD_t^{-1}.
\end{equation}
Observing the relation in \eqref{lower_bounded_matrix} we realize that the lower bound for the eigenvalues of matrix $\bbD_t^{-1}$ in \eqref{lower_bounded_matrix2} holds for the eigenvalues of the Hessian inverse approximation $\hbH_{t}^{-1}$. Therefore, all the eigenvalues of the Hessian inverse approximation $\hbH_{t}^{-1}$ are greater than ${1}/({2(1-\delta)+\alpha M})$. This completes the proof of the claim in \eqref{bounded_Hessian_inverse}.


\section{Proof of Theorem  \ref{linear_convergence}}\label{app_linear_convergecene_theorem}

To prove global convergence of the Network Newton method we first introduce two technical lemmas. In the first lemma we use the result of Lemma \ref{Hessian_Lipschitz_countinous}, namely, that the objective function Hessian $\nabla^2F(\bby)$ is Lipschitz continuous, to develop an upper bound for the objective function value $F(\bby)$ using the first three terms of its Taylor expansion. In the second lemma we construct an upper bound for the objective function error at step $t+1$, namely $F(\bby_{t+1})-F(\bby^{*})$, in terms of the error at step $t$, namely $F(\bby_t)-F(\bby^{*})$.

\begin{lemma}\label{Lemma_taylor_expansion}
Consider the objective function $F(\bby)$ as defined in \eqref{centralized_opt_problem}. If Assumptions 2 and \ref{Lipschitz_assumption} hold true, then for any $\bby,\hby\in\reals^{np}$ the following relation holds
\begin{align}\label{obj_fun_taylor_expan_claim}
F(\hby) & \leq F(\bby) + \nabla F(\bby)^{T}(\hby-\bby)\\
	& \qquad 
		+ \frac{1}{2}(\hby-\bby)^{T}\nabla^{2} F(\bby)(\hby-\bby)
		+ \frac{\alpha L}{6}\|\hby-\bby\|^3.\nonumber
\end{align}

\end{lemma}

\begin{myproof}
Since objective function $F$ is twice differentiable, based on the Fundamental Theorem of Calculus we can write
\begin{equation}\label{fundamental_result_1}
F(\hby)=F(\bby)+\int_{0}^1 \nabla F(\bby+\omega(\hby-\bby))^T (\hby-\bby)\ d\omega,
\end{equation}
where $\nabla F$ is the gradient of function $F$. We proceed by adding and subtracting the term $\nabla F(\bby)^T(\hby-\bby)$ to the right hand side of \eqref{fundamental_result_1} which yields 
\begin{align}\label{taylor_first_two terms}
F(\hby) &=F(\bby)+\nabla F(\bby)^T(\hby-\bby) \\ &\ +\!
	\int_{0}^1\big[ \nabla F(\bby+\omega(\hby-\bby))-\nabla F(\bby)\big]^T (\hby-\bby)\ d\omega. \nonumber 
\end{align}
We apply again the Fundamental Theorem of Calculus but for the gradient $\nabla F$. It follows that for any vectors $\bbz,\hbz \in \reals^{np}$ we can write 
\begin{equation}\label{fundamental_result2}
\nabla F(\hbz)= \nabla F(\bbz)+\int_{0}^1 \nabla^2 F(\bbz+s(\hbz-\bbz)) (\hbz-\bbz)\ ds,
\end{equation}
where $\nabla^2 F$ is the Hessian of function $F$. After setting $\hbz=\bby+\omega(\hby-\bby)$ and $\bbz=\bby$ in \eqref{fundamental_result2} and rearranging terms it follows that
\begin{align}\label{differnce_of_gradients}
&\nabla F(\bby+\omega(\hby-\bby)) -\nabla F(\bby) =  \\ & \int_{0}^1 \nabla^2 F(\bby+s(\bby+\omega(\hby-\bby)-\bby)) (\bby+\omega(\hby-\bby)-\bby)\ ds.\nonumber
\end{align}
Observing the fact that $\bby+\omega(\hby-\bby)-\bby=\omega(\hby-\bby)$, we can further simplify \eqref{differnce_of_gradients} to
\begin{align}\label{differnce_of_gradients_2}
\nabla F(\bby+\omega(\hby-\bby))& -\nabla F(\bby)=\\
&\int_{0}^1 \nabla^2 F(\bby+s\omega(\hby-\bby)) \omega(\hby-\bby)\  ds. \nonumber
\end{align}
Based on the relation for the difference of gradients $\nabla F(\bby+\omega(\hby-\bby)) -\nabla F(\bby)$ in \eqref{differnce_of_gradients_2}, we can rewrite \eqref{taylor_first_two terms} by applying this substitution. This yields
\begin{align}\label{taylor_first_two terms_2}
F(\hby)& =F(\bby)+\nabla F(\bby)^T(\hby-\bby)\\
	&\!\!\!\!\!+ \int_{0}^1\! \int_{0}^1 \! \omega(\hby-\bby)^T \nabla^2 F(\bby+s\omega(\hby-\bby))  (\hby-\bby)\ \! ds\ \!d\omega. \nonumber 
\end{align}
We proceed by adding and subtracting the quadratic integral $\int_{0}^1 \int_{0}^1\omega(\hby-\bby)^T\nabla^2 F(\bby)(\hby-\bby)dsd\omega$ to the right hand side of \eqref{taylor_first_two terms_2} to write
\begin{align}\label{taylor_first_two terms_4}
&F(\hby)\!=\!F(\bby)+\nabla F(\bby)^T(\hby-\bby) \nonumber \\
&\qquad+
	\int_{0}^1 \int_{0}^1 \omega(\hby-\bby)^T\nabla^2 F(\bby)  (\hby-\bby)\ ds \ d\omega \nonumber\\
	&\qquad +  \int_{0}^1 \int_{0}^1 \omega(\hby-\bby)^T\Big[\nabla^2 F(\bby+s\omega(\hby-\bby))-\nabla^2 F(\bby)\Big] \nonumber \\
	&\qquad \qquad \qquad \qquad (\hby-\bby)\ ds  \ d\omega .
\end{align}
Observe that the term $(\hby-\bby)^T\nabla^2 F(\bby)  (\hby-\bby)$ in the third summand of \eqref{taylor_first_two terms_4} is not a function of $\omega$ or $s$. Hence, we can move this term outside of the integral and simplify the integral to $\int_{0}^1 \int_{0}^1 \omega \ dsd\omega=1/2$. As a result of these observation the third summand of \eqref{taylor_first_two terms_4} can be replaced by $(1/2)(\hby-\bby)^T\nabla^2 F(\bby)  (\hby-\bby)$ and we can rewrite \eqref{taylor_first_two terms_4} as
\begin{align}\label{taylor_first_three_terms}
F(\hby)&=F(\bby)\!+\!\nabla F(\bby)^T(\hby-\bby)\!+\!
	\frac{1}{2}(\hby-\bby)^T\nabla^2 F(\bby)  (\hby-\bby) \nonumber\\
	 &+  \int_{0}^1 \int_{0}^1 \omega(\hby-\bby)^T\Big[\nabla^2 F(\bby+s\omega(\hby-\bby))-\nabla^2 F(\bby)\Big] \nonumber \\
	&\qquad \qquad \qquad \qquad (\hby-\bby)\ ds  \ d\omega .
\end{align}
We proceed now to construct an upper bound for the integral in \eqref{taylor_first_three_terms}. Observe that according to the definition of the Euclidean norm of a matrix we have the inequality $(\hby-\bby)^T\big[\nabla^2 F(\bby+s\omega(\hby-\bby))-\nabla^2 F(\bby)\big]  (\hby-\bby) \leq \|\nabla^2 F(\bby+s\omega(\hby-\bby))-\nabla^2 F(\bby) \|_2 \| \hby-\bby\|^2$. By applying this inequality we have an upper bound for the integral in \eqref{taylor_first_three_terms} that results in
\begin{align}\label{taylor_first_three_terms_2}
&F(\hby)\!\leq\! F(\bby)\!+\!\nabla F(\bby)^T(\hby-\bby)\!+\!
	\frac{1}{2}(\hby-\bby)^T\nabla^2 F(\bby)  (\hby-\bby) \nonumber\\
	& + \| \hby-\bby\|^2\!\! \int_{0}^1\!\! \int_{0}^1 \! \omega\|\nabla^2 F(\bby+s\omega(\hby-\bby))\!-\!\nabla^2 F(\bby) \|_2 dsd\omega.
\end{align}
The next step is to provide an upper bound for the term $\|\nabla^2 F(\bby+s\omega(\hby-\bby))\!-\!\nabla^2 F(\bby) \|_2$ in the right hand side of \eqref{taylor_first_three_terms_2}. Lemma \ref{Hessian_Lipschitz_countinous} shows that the penalized objective function Hessian $\nabla^2F$ is Lipschitz continuous with parameter $\alpha L$. Therefore, we can write 
\begin{align}\label{upper_bound_for_Hessian_difference}
\!\!\left\|  \nabla^2 F(\bby\!+\!s\omega(\hby\!-\!\bby))\!-\!\nabla^2 F(\bby)  \right\|_{2} & 
		\!\leq \! \alpha L \| \bby\!+\!s\omega(\hby-\bby)\!-\!\bby \|  \nonumber \\
	&   = \alpha L s \omega\|\hby-\bby  \|           .
\end{align} 
By considering \eqref{taylor_first_three_terms_2} and substituting $\|  \nabla^2 F(\bby+s\omega(\hby-\bby))-\nabla^2 F(\bby)  \|_{2} $ by the upper bound in \eqref{upper_bound_for_Hessian_difference}, we obtain 
\begin{align}\label{taylor_first_three_terms_3}
F(\hby) &\leq F(\bby)+\nabla F(\bby)^T(\hby-\bby)\!+\!
	\frac{1}{2}(\hby-\bby)^T\nabla^2 F(\bby)  (\hby-\bby)  \nonumber \\
  & \qquad + \| \hby-\bby\|^2 \int_{0}^1 \int_{0}^1 \alpha L s \omega^2\|\hby-\bby  \|   \ dsd\omega.
\end{align}
Now observe that since $\|\hby-\bby  \|$ does not depend on $s$ or $\omega$, the integral in the last summand of \eqref{taylor_first_three_terms_3} can be simplified as
\begin{align}\label{simplification_of_integral}
 \int_{0}^1\!\! \int_{0}^1 \alpha L s \omega^2\|\hby-\bby  \|   \ dsd\omega 
		&= \alpha L\| \hby-\bby\|  \int_{0}^1\!\! \int_{0}^1  s \omega^2  \ dsd\omega \nonumber\\
		&= \frac{\alpha L}{6}\| \hby-\bby\|.
\end{align} 
The simplification in \eqref{simplification_of_integral} for the last summand of \eqref{taylor_first_three_terms_3} implies the claim in \eqref{obj_fun_taylor_expan_claim} is valid.  
\end{myproof}


Lemma \ref{Lemma_taylor_expansion} shows an upper bound for the Taylor expansion of the objective function value $F(\hby)$. We use the result of Lemma \ref{Lemma_taylor_expansion} to establish an upper bound for the objective function error at step $t+1$ in terms of the error at step $t$. This result is proven in the following lemma.

\begin{lemma}\label{important_lemma_for_convergnce}
Consider the NN-$K$ method as defined in \eqref{diagonal_matrix}-\eqref{update_formula_NN} and the objective function $F(\bby)$ as defined in \eqref{centralized_opt_problem}. Further, recall the definition of $\bby^{*}$ as the optimal argument of the objective function $F(\bby)$. If assumptions \ref{ass_weight_bounds}, \ref{convexity_assumption}, and \ref{Lipschitz_assumption} hold true, the sequence of objective function value errors $\{F(\bby_{t}) - F(\bby^*)\}$ satisfies 
\begin{align}\label{claim_new}
\!F(\bby_{t+1}) -F(\bby^*)  
		&\leq \left[ 1 -{ \left(2\eps-{\epsilon}^2\right) \alpha m\lambda}  \right] \! [F(\bby_{t}) - F(\bby^*)]  \nonumber \\
	& \qquad +  \frac{\alpha L\epsilon^3\Lambda^{3}}{6\lambda^{\frac{3}{2}}} {\left[F(\bby_{t}) -F(\bby^*) \right]}^{\frac{3}{2}}.
\end{align}
\end{lemma}


\begin{myproof}
Recall the result of Lemma \ref{Lemma_taylor_expansion}. By setting $\hby:=\bby_{t+1}$ and $\bby:=\bby_t$ in \eqref{obj_fun_taylor_expan_claim} we obtain
\begin{align}\label{obj_fun_taylor_expan2}
F(\bby_{t+1}) & \leq F(\bby_{t}) + \bbg_t^{T}(\bby_{t+1}-\bby_{t})  \\
	&\quad+\! \frac{1}{2}(\bby_{t+1}\!-\!\bby_{t})^{T}\bbH_t(\bby_{t+1}\!-\!\bby_{t})  
	\!+\! \frac{\alpha L}{6}\|\bby_{t+1}\!-\!\bby_{t}\|^3\!, \nonumber 
\end{align}
where $\bbg_t:=\nabla F(\bby_{t}) $ and $\bbH_t:=\nabla^{2} F(\bby_t)$. From the definition of the NN-$K$ update formula in \eqref{Hessian_approximation_iteration} we can write the difference of two consecutive variables as $\bby_{t+1}-\bby_t=-\eps\hbH_{t}^{-1} \bbg_t$. Making this substitution in \eqref{obj_fun_taylor_expan2} implies
\begin{align}\label{obj_fun_taylor_expan3}
F(\bby_{t+1}) &\leq F(\bby_{t})
	 - \epsilon \bbg_t^{T}\hbH_{t}^{-1}\bbg_t
	 + \frac{\epsilon^2}{2} \bbg_t^{T}\hbH_{t}^{-1}\bbH_{t}							\hbH_{t}^{-1} \bbg_t
	 \nonumber \\
&\qquad	
	+ \frac{\alpha L\epsilon^3}{6}\|\hbH_{t}^{-1}\bbg_t\|^3.
\end{align}
According to the definition of error matrix $\bbE_t$ in \eqref{error_matrix}, we can substitute $\hbH_{t}^{-1/2}\bbH_t\hbH_{t}^{-1/2}$ by $\bbI-\bbE_t$. By making this substitution into the third summand of \eqref{obj_fun_taylor_expan3} we obtain
\begin{align}\label{obj_fun_taylor_expan4}
F(\bby_{t+1}) &\leq F(\bby_{t})
	 - \epsilon \bbg_t^{T}\hbH_{t}^{-1}\bbg_t
	\!+\! \frac{\epsilon^2}{2}\bbg_t^{T}\hbH_{t}^{-\frac{1}{2}}							(\bbI-\bbE_t)\hbH_{t}^{-\frac{1}{2}} \bbg_t\nonumber \\
	& \qquad
	+ \frac{\alpha L\epsilon^3}{6}\|\hbH_{t}^{-1} \bbg_t\|^3.
\end{align}
Proposition \ref{error_matrix_proposition} shows that the error matrix $\bbE_t $ is always positive semidefinite. As a result, we conclude that the quadratic form $ \bbg_t^{T}\hbH_{t}^{-{1}/{2}}\bbE_t\hbH_{t}^{-{1}/{2}}  \bbg_t$ is always nonnegative. Considering this lower bound we can simplify \eqref{obj_fun_taylor_expan4} to 
\begin{equation}\label{obj_fun_taylor_expan5}
F(\bby_{t+1}) \leq F(\bby_{t})
	 - \frac{\left(2\epsilon-\epsilon^{2}\right) }{2} \bbg_t^{T}\hbH_{t}^{-1} \bbg_t
	+ \frac{\alpha L\epsilon^3}{6}\|\hbH_{t}^{-1}  \bbg_t\|^3.
\end{equation}
Note that since the stepsize is not larger than $1$, we obtain that $2\eps - \eps^{2}$ is positive. Moreover, recall the result of Lemma \ref{Hessian_inverse_eigenvalue_bounds_lemma} that all the eigenvalues of the Hessian inverse approximation $\hbH_{t}^{-1}$ are lower and upper bounded by $\lambda$ and $\Lambda$, respectively. These two observations imply that we can replace the term $ \bbg_t^{T}\hbH_{t}^{-1} \bbg_t$ by its lower bound $ \lambda\|\bbg_t\|^2$. Moreover, existence of upper bound $\Lambda$ for the eigenvalues of Hessian inverse approximation $\hbH_{t}^{-1}$ implies that the term $\|\hbH_{t}^{-1}  \bbg_t\|^3$ is upper bounded by $\Lambda^3\|\bbg_{t}\|^3$. Substituting these bounds for the second and third terms of \eqref{obj_fun_taylor_expan5} and subtracting the optimal objective function value $F(\bby^*)$ from both sides of inequality \eqref{obj_fun_taylor_expan5} leads to
\begin{align}\label{obj_fun_taylor_expan6}
F(\bby_{t+1})-F(\bby^*)  & \leq F(\bby_{t})-F(\bby^*)
	- \frac{\left(2\epsilon-\epsilon^{2}\right)\lambda}{2}\|\bbg_t\|^{2}  \nonumber \\
	& \qquad + \frac{\alpha L\epsilon^3\Lambda^{3}}{6}\|\bbg_t\|^3.
\end{align}
We now find lower and upper bounds for the norm of gradient $\|\bbg_t\|$  in terms of the objective function error $F(\bby_t)-F(\bby^*)$. As it follows from Proposition \ref{eigenvalue_bounds}, the eigenvalues of Hessian $\bbH_t$ are bounded by $\alpha m$ and $2+\alpha M$. Taking a Taylor expansion of the objective function $F(\bby)$ around $\bbw$ and using the lower bound $\alpha m $ for the Hessian eigenvalues yields
\begin{equation}\label{taylor_lower_bound}
   F(\bby) \geq\ F(\bbw) +\nabla F(\bbw)^{T}(\bby-\bbw)   
   + {{\alpha m}\over{2}}\|{\bby - \bbw}\|^{2}.
\end{equation}
For fixed $\bbw$, the right hand side of \eqref{taylor_lower_bound} is a quadratic function of $\bby$ whose minimum argument we can find by setting its gradient to zero. Doing this yields the minimizing argument $\hby = \bbw- (1/m) \nabla  F(\bbw)$ implying that for all $\bby$ we must have
\begin{alignat}{2}\label{lower_bound_for_gradient}
F(\bby) \geq\ 
    &\ F(\bbw) +\nabla F(\bbw)^{T}(\hby-\bbw)   
   + {{\alpha m}\over{2}}\|{\hby - \bbw}\|^{2} \nonumber \\
   \ =\ 
         &\ F(\bbw) - \frac{1}{2\alpha m} \| \nabla F(\bbw)\|^{2} .
\end{alignat}
The bound in \eqref{lower_bound_for_gradient} is true for all $\bbw$ and $\bby$. In particular, for $\bby=\bby^{*}$ and $\bbw=\bby_{t}$ \eqref{lower_bound_for_gradient} yields
\begin{equation}\label{transition}
   F(\bby^*) \geq  F(\bby_t) - \frac{1}{2\alpha m} \| \nabla F(\bby_t)\|^{2}.
\end{equation} 
Rearrange terms in \eqref{transition} to obtain $2\alpha m (F(\bby_t) -F(\bby^*))$ as a lower bound for $\| \nabla F(\bby_t)\|^{2}=\|\bbg_t\|^2$. Now substitute the lower bound $2\alpha m (F(\bby_t) -F(\bby^*))$ for squared norm of gradient $\| \bbg_t\|^{2}$ in the second summand of \eqref{obj_fun_taylor_expan6} to obtain 
\begin{align}\label{obj_fun_taylor_expan7}
	F(\bby_{t+1})\!-\!F(\bby^*)  & \!\leq \left[1
	\! -\! {\left(2\epsilon-\epsilon^{2}\right)\alpha m\lambda}\right]\! (F(\bby_{t})\!-\!F(\bby^*)) \nonumber \\
	& \qquad + \frac{\alpha L\epsilon^3\Lambda^{3}}{6}\|\bbg_t\|^3.
\end{align}
To find an upper bound for the norm of gradient $\|\bbg_t\|$ in terms of the objective function error $F(\bby_t )-F(\bby^*) $ we first use the Taylor expansion of the objective function $F(\bby)$ around $\bbw$ by considering the fact that $2(1-\delta)+\alpha M$ is an upper bound for the eigenvalues of the Hessian. For any vectors $\hby$ and $\bby$ in $\reals^{np}$ we can write
\begin{equation}\label{taylor_upper_bound}
   F(\bby) \leq F(\hby) +\nabla F(\hby)^{T}(\bby-\hby)   
   + {{{2(1-\delta)+\alpha M}\over{2}}\|{\bby - \hby}\|^{2}.}
\end{equation}
Notice that according to the definition of $\lambda $ in \eqref{definition_of_lambdas} we can substitute $2(1-\delta)+\alpha M$ by $1/\lambda$. Implementing this substitution and minimizing both sides of the equality with respect to $\bby$ yields
\begin{equation}\label{taylor_upper_bound_2}
   F(\bby^*) \leq\ F(\hby) -\lambda \| \nabla F(\hby) \|^2.
\end{equation}
Setting $\hby=\bby_{t}$, observing that by definition  $ \| \nabla F(\bby_{t}) \|=\|\bbg_t\|$, rearranging terms, and taking the square root of both sides of the resulting inequality leads to 
\begin{equation}\label{taylor_upper_bound_3}
 \| \bbg_t\| \leq  \left[\frac{1}{\lambda}  \left[F(\bby_t) - F(\bby^*)\right] \right]^{\frac{1}{2}} .
\end{equation}
Replacing the upper bound in \eqref{taylor_upper_bound_3} for the norm of the gradient $ \| \bbg_t\| $ in the last term of \eqref{obj_fun_taylor_expan7} yields the claim in \eqref{claim_new}.
\end{myproof}

We use the result of Lemma \ref{important_lemma_for_convergnce} to prove linear convergence of the sequence of objective function errors $F(\bby_t)-F(\bby^*)$ to zero.\\

\textbf{Proof of Theorem \ref{linear_convergence}:} 
To simplify upcoming derivations define the sequence $\beta_t$ as
\begin{equation}\label{beta_definition}
\beta_{t} \!:=\!    {(2-{\epsilon})\epsilon \alpha m\lambda} -  \frac{\epsilon^3\alpha L\Lambda^{3}\left[{F(\bby_{t})\!-\!F(\bby^{*}) }\right]^{\frac{1}{2}}}{6\lambda^{\frac{3}{2}}} .
\end{equation}
Recall the result of Lemma \ref{important_lemma_for_convergnce}. Factorizing  $F(\bby_{t})-F(\bby^{*})$ from the terms of the right hand side of \eqref{claim_new} in association with the definition of $\beta_{t}$ in \eqref{beta_definition} implies that we can simplify \eqref{claim_new} as
\begin{equation}\label{obj_fun_taylor_expan100}
F(\bby_{t+1}) -F(\bby^*) \leq (1-\beta_{t})  {\left(F(\bby_{t}) -F(\bby^*) \right)}.
\end{equation}
To prove global convergence of objective function error $F(\bby_{t}) -F(\bby^*)$ we need to show that for all time steps $t$, the constants $\beta_{t}$ are strictly smaller than $1$ and larger than $0$, i.e., that $0<\beta_t<1$ for all times $t$. 

We first show that $\beta_{t}$ is less than $1$ for all $t\geq0$. To do so observe that the second term in the right and side of \eqref{beta_definition} is nonnegative. It is therefore true that 
\begin{equation}\label{beta_definition_bound}
   \beta_{t}  \leq   {(2-{\epsilon})\epsilon \alpha m\lambda}. 
\end{equation}
Considering the inequality $(\epsilon-1)^2\geq0$ it is trivial to derive that $\epsilon(2-\epsilon)\leq1$. Moreover, considering the facts that $m<M$ and $1-\delta>0$, we obtain $\alpha m < \alpha M +(1-\delta )$ which yields $\alpha m / (\alpha M +2(1-\delta ))<1$. Considering the definition of $\lambda $ in \eqref{definition_of_lambdas} we can substitute $1/(2(1-\delta)+\alpha M)$ by $\lambda$ and write $\alpha m\lambda<1$.
By multiplying these two ratios, both of which are smaller than $1$, we conclude that
\begin{equation}\label{eps_cond6}
(2-\epsilon)\epsilon\alpha m\lambda <1.
\end{equation}
That $\beta_t<1$ follows by combining \eqref{beta_definition_bound} with \eqref{eps_cond6}.

To prove that $0<\beta_t$ for all $t\geq0$ we prove that this is true for $t=0$ and then prove that the $\beta_t$ sequence is increasing. To show that $\beta_0$ is positive first note that since the stepsize $\epsilon$ satisfies the condition in \eqref{step_size_condition} we can write  
\begin{equation}\label{eps_cond1}
\epsilon \leq \left[\frac{3m\lambda^{\frac{5}{2}}}{ L\Lambda^{3}{{(F(\bby_{0})-F(\bby^{*}) )^{\frac{1}{2}}}} } \right]^{\frac{1}{2}},
\end{equation}
By computing the squares of both sides of \eqref{eps_cond1}, multiplying the right hand side of the resulting inequality by 2 to make the inequality strict, and factorizing ${\alpha m}\lambda$ from the term in the resulting right hand side we obtain
\begin{equation}\label{eps_cond2}
\epsilon^2  < \frac{6\lambda^{\frac{3}{2}}}{\alpha L\Lambda^{3}[{F(\bby_{0})-F(\bby^{*}) }]^{\frac{1}{2}}} \times {\alpha m\lambda}.
\end{equation}
If we now divide both sides of the inequality in \eqref{eps_cond2} by the first multiplicand in the right hand side of \eqref{eps_cond2} we obtain
\begin{equation}\label{eps_cond3}
 \frac{\epsilon^2\alpha L\Lambda^{3}[{F(\bby_{0})-F(\bby^{*}) }]^{\frac{1}{2}}}{6\lambda^{\frac{3}{2}}} < {\alpha m\lambda}.
 \end{equation}
Observe that based on the hypothesis in \eqref{step_size_condition} the step size $\epsilon$ is smaller than $1$ and it is then trivially true that $2-\epsilon\geq1$. This observation shows that if we multiply the right hand side of \eqref{eps_cond3}  by $2(1-\epsilon/2)$ the inequality still holds,
\begin{equation}\label{eps_cond4}
 \frac{\epsilon^2\alpha L\Lambda^{3}{({F(\bby_{0})-F(\bby^{*}))^{\frac{1}{2}}} }}{6\lambda^{\frac{3}{2}}} 
 <  {\alpha m  (2-{\epsilon})\lambda}.
\end{equation}
Furhter multiplying both sides of inequality \eqref{eps_cond4} by $\epsilon$ and rearranging terms leads to
\begin{equation}\label{eps_cond5}
 {\alpha m \eps (2-{\epsilon})\lambda} -  \frac{\epsilon^3\alpha L \Lambda^3[{F(\bby_{0})-F(\bby^{*}) }]^{\frac{1}{2}}}{6\lambda^{\frac{3}{2}}} \!>\!0.
\end{equation}
According to the definition of $\beta_{t}$ in \eqref{beta_definition}, the result in \eqref{eps_cond5} implies that $\beta_{0}>0$. 

Observing that $\beta_0$ is positive, to show that for all $t$ the sequence of $\beta_{t}$ is positive it is sufficient to prove that the sequence $\beta_t$ is increasing, i.e., that $\beta_{t}<\beta_{t+1}$ for all $t$. We use strong induction to prove $\beta_{t}<\beta_{t+1}$ for all $t\geq0$. By setting $t=0$ in \eqref{obj_fun_taylor_expan100} the inequality can be written as 
\begin{equation}\label{inequality_for_t=0}
F(\bby_{1}) -F(\bby^*)\leq (1-\beta_{0})(F(\bby_{0}) -F(\bby^*)).
\end{equation}
Considering the result in \eqref{inequality_for_t=0} and the fact that $0<\beta_{0}<1$, we obtain that the objective function error at time $t=1$ is strictly smaller than the error at time $t=0$, i.e.
\begin{equation}\label{inequality_for_t=0_2}
F(\bby_{1}) -F(\bby^*)<F(\bby_{0}) -F(\bby^*).
\end{equation}
Observe now that in the definition of sequence $\beta_{t}$ in \eqref{beta_definition} the objective function error term $F(\bby_t) -F(\bby^*)$ appears in the numerator of negative term. Therefore, a smaller objective function error $F(\bby_t) -F(\bby^*)$ leads to a larger coefficient $\beta_t$. Hence, this observation in association with the result in \eqref{inequality_for_t=0_2} leads to the conclusion, 
\begin{equation}\label{khoob_tarin_result}
\beta_0<\beta_1.
\end{equation}
To complete the strong induction argument assume now that $\beta_{0}<\beta_{1}<\dots<\beta_{t-1}<\beta_{t}$ and proceed to prove that if this is true we must have $\beta_{t}<\beta_{t+1}$. Begin by observing that since $0<\beta_{0}$ the induction hypothesis implies that for all $u\in\{0,\dots,t\}$ the constant $\beta_{u}$ is also positive, i.e., $0<\beta_{u}$. Further recall that for all $t$ the sequence $\beta_{t}$ is also smaller than $1$ as already proved. Combining these two observations we can conclude that $0<\beta_{u}<1$ for all $u\in\{0,\dots,t\}$. Consider now the inequality in \eqref{obj_fun_taylor_expan100} and utilize the fact that $0<\beta_{u}<1$ for all $u\in\{0,\dots,t\}$ to conclude that
\begin{equation}\label{obj_fun_taylor_expan1000}
F(\bby_{u+1}) -F(\bby^*)<  F(\bby_{u}) -F(\bby^*),
\end{equation}
for all $u\in\{0,\dots,t\}$. Setting $u=t$ in \eqref{obj_fun_taylor_expan1000} we conclude that $F(\bby_{t+1}) -F(\bby^*)<F(\bby_{t}) -F(\bby^*)$. By further repeating the argument leading from \eqref{khoob_tarin_result} to \eqref{inequality_for_t=0_2} we can conclude that 
\begin{equation}\label{khoob_tarin_result_2}
   \beta_t<\beta_{t+1}.
\end{equation}
The strong induction proof is complete and we can now claim that for all times $t$ 
\begin{equation}\label{eqn_beta_chain}
   0<\beta_{0}<\beta_{1}<\dots<\beta_{t}<1.
\end{equation}
The relationship in \eqref{obj_fun_taylor_expan100} and the property in \eqref{eqn_beta_chain} imply convergence of the objective function value sequence to the optimal argument, i.e.  $\lim_{t\to \infty} F(\bby_{t}) -F(\bby^*)=0$. To conclude that the convergence rate is at least linear simply observe that if the sequence $\beta_{t}$ is increasing as per \eqref{eqn_beta_chain}, the sequence $1-\beta_t$ is decreasing and satisfies 
\begin{equation}\label{inter_milan_champion}
   0<1-\beta_{t}< 1-\beta_{0}<1,
\end{equation}
for all time steps $t$. Applying the inequality in \eqref{obj_fun_taylor_expan100} recursively and considering the inequality in \eqref{inter_milan_champion} yields 
\begin{equation}\label{obj_fun_taylor_expan200}
   F(\bby_{t}) -F(\bby^*) \leq (1-\beta_{0})^t  {\left(F(\bby_{0}) -F(\bby^*) \right)},
\end{equation}
which shows the objective function error sequence $F(\bby_{t}) -F(\bby^*)$ converges to $0$ at least linearly with constant $(1-\beta_{0})$. By setting $\zeta=\beta_{0}$, the claim in \eqref{linear_convegrence_claim} follows.

\end{appendices}
\bibliographystyle{IEEEtran}
  \bibliography{bmc_article}
   \end{document}